\documentclass[12pt]{amsart}
\usepackage{euscript,amssymb,amscd,mathrsfs}
\usepackage[matrix,arrow,curve]{xy}
\makeatletter
\@addtoreset{equation}{subsection}
\makeatother

\renewcommand{\theequation}{\thesubsection.\arabic{equation}}

\newcommand{\para}
{\refstepcounter{subsubsection}\refstepcounter{equation}\par\smallskip\noindent
{\rm \bf (\theequation)\hspace{0.2cm}}}

\newcommand{\comment}[1]{}
\newcommand{\muu}{{\boldsymbol\mu}}
\newcommand{\ep}{\varepsilon}
\newcommand{\var}{\varphi}
\newcommand{\down}[1]{\left\lfloor #1\right\rfloor}
\newcommand{\fr}[1]{\left\{ #1 \right\}}
\newcommand{\up}[1]{\left\lceil #1\right\rceil}
\newcommand{\RR}{{\mathbb R}}
\newcommand{\CC}{{\mathbb C}}
\newcommand{\QQ}{{\mathbb Q}}
\newcommand{\R}{{\mathfrak{R}}}
\newcommand{\PP}{{\mathbb P}}

\newcommand{\PPP}{{\mathscr{P}}}
\newcommand{\NNN}{{\mathscr{N}}}
\newcommand{\MMM}{{\mathscr{M}}}
\newcommand{\BBB}{{\mathscr{B}}}

\newcommand{\h}{\mathrm{h}}
\newcommand{\ver}{\mathrm{v}}
\newcommand{\ov}[1]{\overline{#1}}
\newcommand{\comp}{\mathrel{\scriptstyle{\circ}}}

\newcommand{\bD}{\mathbf D}
\newcommand{\bK}{\mathbf K}

\newcommand{\red}{\operatorname{red}}

\newcommand{\codim}{\operatorname{codim}}
\newcommand{\Center}{\operatorname{Center}}

\newcommand{\lcm}{\operatorname{lcm}}

\newcommand{\Const}{\operatorname{Const}}
\newcommand{\Supp}{\operatorname{Supp}}
\newcommand{\Sing}{\operatorname{Sing}}
\newcommand{\smod}{_{\operatorname{mod}}}
\newcommand{\sdiv}{_{\operatorname{div}}{}}
 \newcommand{\ZZ}{\mathbb Z}

 \newcommand{\Diff}{\operatorname{Diff}}

\newcommand{\discr}{\operatorname{discr}}
\newcommand{\totaldiscr}{\operatorname{totaldiscr}}
\newcommand{\Pic}{\operatorname{Pic}}
\newcommand{\Div}{\operatorname{Div}}
\newcommand{\NE}{{\overline{\operatorname{NE}}}}

\newcommand{\q}[1]{\mathbin{\sim_{\scriptscriptstyle{#1}}}}

\renewcommand{\thefootnote}{\fnsymbol{footnote}}
\newcommand{\xref}[1]{{\rm \ref{#1}}}

\newtheorem{theorem}[subsection]{Theorem}
\newtheorem{conjecture}[subsection]{Conjecture}
\newtheorem{proposition}[subsection]{Proposition}
\newtheorem{lemma}[subsection]{Lemma}
\newtheorem{lemma-definition}[subsection]{Lemma-Definition}
\newtheorem{claim}[subsection]{Claim}
\newtheorem{claim1}[subsubsection]{Claim}
\newtheorem{corollary}[subsection]{Corollary}

\theoremstyle{definition}
\newtheorem{definition}[subsection]{Definition}
\newtheorem{assumption}[subsection]{Assumption}
\newtheorem{example}[subsection]{Example}

\newtheorem{remark}[subsection]{Remark}
\newtheorem{remark1}{Remark}[subsection]
\newtheorem*{remark*}{Remark}
\newtheorem{edefinition}[subsection]{}

\newcommand{\com}[1]{}
\newtheorem{Addendum}[subsection]{Addendum}

\author{Yu.~G.~Prokhorov}
\thanks {The first author was 
partially supported by grants CRDF-RUM, No. 1-2692-MO-05 and
RFBR, No. 05-01-00353-a, 06-01-72017.
The second author was
partially supported by NSF grant DMS-0400832.}

\address{Yu. G. Prokhorov:
Department of Algebra, Faculty of Mathematics, Moscow State
University, Moscow 117234, Russia}
\email{prokhoro@mech.math.msu.su}

\author{V.~V.~Shokurov}
\address{V.~V.~Shokurov: The Johns Hopkins University,
Department of Mathematics, Baltimore, Maryland 21218, USA
\phantom{ggggggggggg}
\linebreak
\phantom{\qquad}
Steklov Mathematical Institute, Russian Academy of Sciences,
Gubkina str. 8, 119991, Moscow, Russia}
\email{shokurov@math.jhu.edu}

\date{}

\title[The second theorem on complements]{Towards the second main theorem on complements}

\begin{document}
\maketitle 
\begin{abstract}
We prove the boundedness 
of complements modulo two conjectures:
Borisov-Alexeev conjecture and effective adjunction
for fibre spaces.
We discuss the last conjecture and prove it in two
particular cases. 
\end{abstract}
\tableofcontents
\section{Introduction}
This paper completes our previous work \cite{Prokhorov-Shokurov-2001} modulo two
conjectures \xref{BAB} and \xref{conj-main-adj}. The first one relates to
Alexeev's, A. and L. Borisov's conjecture:

\begin{conjecture}
\label{BAB}
Fix a real number $\ep>0$. Let $(X,B=\sum b_iB_i)$ be a $d$-dimensional log
canonical pair with nef $-(K_X+B)$, that is, $(X,B)$ is a log semi-Fano
variety \textup(cf. Definition \xref{definition}\textup). Assume also that 
\begin{enumerate}
\item
$K+B$ is $\ep$-lt\textup; and
\item
$X$ is FT that is $(X,\Theta)$ is a klt log Fano variety
with respect to some boundary $\Theta$.
\end{enumerate}
\par\medskip\noindent
Then $X$ is bounded in the moduli sense, i.e., it belongs to an
algebraic family $\mathcal{X}(\ep,d)$.
\end{conjecture}
This conjecture was proved in dimension $2$ by V. Alexeev
\cite{Alexeev-1994} and in toric case by A. Borisov and L. Borisov
\cite{Borisov-Borisov} (see also \cite{Nikulin-1990-e}, \cite{Borisov-1996},
\cite{Borisov-2001}, \cite{McKernan-2002-AG}).

\begin{remark}
We hope that Conjecture~\ref{BAB} can be generalized by weakening
condition (ii). For example, we hope that instead of (ii) one can assume that
\begin{enumerate}
\item[(ii)${}'$]
$X$ is rationally connected,
cf. \cite{Zhang-Qi-2006}.
\end{enumerate}
\end{remark}

Recall that a log pair $(X,B)$ is said to be \emph{$\ep$-log
terminal} (or simply \emph{$\ep$-lt}) if
$\totaldiscr (X,B)>-1+\ep$, see Definition \ref{def-ep-lt} below.

\begin{theorem}[{\cite{Alexeev-1994}}]
\label{BAB_th}
Conjecture \xref{BAB} holds in dimension two.
\end{theorem}

The second conjecture concerns with Adjunction Formula and will be
discussed in Section~\xref{sect-Adj}.

Our main result is the
following.

\begin{theorem}
\label{main-result}
Fix a finite subset $\R\subset [0,\, 1]\cap \QQ$.
Let $(X,B)$ be a klt
log semi-Fano variety of dimension $d$ such that
$X$ is FT and
the multiplicities of $B$ are contained in $\Phi(\R)$
\textup(see \xref{not-hyperstand}\textup).
Assume the LMMP in dimension $\le d$.
Further, assume that Conjectures~\xref{BAB} and \xref{conj-main-adj}
hold in dimension $\le d$.
Then $K+B$ has bounded complements. More precisely,
there is a positive integer $n=n(d,\R)$ divisible by 
denominators of all $r\in \R$ and such that
$K+B$ is
$n$-complemented. Moreover, $K+B$ is $nI$-complemented
for any positive integer $I$.

In particular, 
\[
\left|-nK-nS -\down{(n+1)D}\right|\neq \emptyset,
\]
where $S:=\down B$ and $D:=B-S$.
For $B=0$, 
$|-nK|\neq \emptyset$, where
$n$ depends only on $d$.
\end{theorem}
Note that the last paragraph is an immediate consequence
of the first statemet and the definition of complements.

In the case when $K+B$ is numerically trivial our result
is stronger. For the definition of $0$-pairs we refer to
\xref{definition}.

\begin{theorem}[cf. {\cite{Blache-1995}}, {\cite{Ishii-2000-AG}}]
\label{main-result0}
Fix a finite subset $\R\subset [0,\, 1]\cap \QQ$.
Let $(X,B)$ be a $0$-pair of dimension $d$ such that
$X$ is FT and the multiplicities of $B$ are contained in $\Phi(\R)$.
Assume the LMMP in dimension $\le d$.
Further, assume that Conjectures~\xref{BAB} and \xref{conj-main-adj}
hold in dimension $\le d$.
Then there is a positive integer
$n=n(d,\R)$ such that $n(K_X+B)\sim 0$.
\end{theorem}

\begin{Addendum}
\label{addendum-1}
Our proofs show that we do not need Conjecture~\xref{BAB} in
dimension $d$ in all
generality. We need it only for some special value 
$\ep^o:=\ep^o(d,\R)>0$
\textup(cf. Corollary \xref{main-result-3-Cor-1-b}\textup).
\end{Addendum}

We will prove Conjecture 
\xref{conj-main-adj} in \S \ref{section-Adj-conj-part-cases}
in dimension $\le 3$
under additional assumption that the total space is projective
and FT.

\begin{corollary}
\label{main-result-3-Cor-1-b}
Fix a finite rational subset $\R\subset [0,\, 1]\cap \QQ$
and let $I$ be the least common multiple of denominators of 
all $r\in \R$.
Let $(X,B)$ be a klt log semi-Fano threefold such that
$X$ is FT and
the multiplicities of $B$ are contained in $\Phi(\R)$.
Assume that Conjecture~\xref{BAB} 
holds in dimension $3$ for $\ep^o$ as in Addendum \xref{addendum-1}.
Then $K+B$ has a bounded $n$-complement such that $I\mid n$.
In particular, there exists a positive integer number $n$ such that 
$I\mid n$ and
\[
\left|-nK-nS -\down{(n+1)D}\right|\neq \emptyset,
\]
where $S:=\down B$ and $D:=B-S$\textup;
$n$ depends only on $\R$ and $\ep^o$ for Addendum \xref{addendum-1}.
For $B=0$, 
$|-nK|\neq \emptyset$, where
$n$ depends only on $\ep^o$.
\end{corollary}
\begin{proof} 
Immediate by Addendum \xref{addendum-1}, Theorem \xref{main-result}, and 
Corollary \xref{corollary-815}.
\end{proof}

\begin{corollary}[cf. {\cite{Shokurov-2000}}]
\label{main-result-2}
Fix a finite subset $\R\subset [0,\, 1]\cap \QQ$
and let $I$ be the least common multiple of denominators of 
all $r\in \R$.
Let $(X,B)$ be a klt log semi-del Pezzo surface
such that $X$ is FT and
the multiplicities of $B$ are contained in $\Phi(\R)$.
Then $K+B$ has a bounded $n$-complement such that $I\mid n$.
In particular, there exists a positive integer $n$ such that $I\mid n$ and 
\[
\left|-nK-nS -\down{(n+1)D}\right|\neq \emptyset,
\]
where $S:=\down B$ and $D:=B-S$\textup;
$n$ depends only on $\R$.
For $B=0$, 
$|-nK|\neq \emptyset$, where
$n$ is an absolute constant.
\end{corollary}
\begin{proof} 
Immediate by Theorems \xref{main-result}, \xref{BAB_th}, and \xref{th-n-n-1}.
\end{proof}

The following corollaries are consequences of 
our techniques. 
The proofs will be given in \ref{sketch-proof-main-result-3-Cor-1-aC},
\ref{sketch-proof-main-result-3-Cor-1-a},
and \ref{sketch-proof-main-result-3-Cor-2div}.

\begin{corollary}[cf. {\cite{Ishii-2000-AG}}]
\label{main-result-3-Cor-1-a}
Fix a finite rational subset $\R\subset [0,\, 1]\cap \QQ$.
Let $(X,D)$ be a three-dimensional $0$-pair such that
$X$ is FT and the multiplicities of $B$ are contained in $\Phi(\R)$.
Assume that $(X,D)$ is not klt.
Then there exists a positive integer $n$ such that
$n(K+D)\sim 0$\textup; this
$n$ depends only on $\R$.
\end{corollary}

\begin{corollary}
\label{main-result-3-Cor-1-2div}
Fix a finite rational subset $\R\subset [0,\, 1]\cap \QQ$
and let $I$ be the least common multiple of denominators of 
all $r\in \R$.
Let $(X,B)$ be a klt log semi-Fano threefold such that
$X$ is FT and
the multiplicities of $B$ are contained in $\Phi(\R)$.
Then there exists a real number $\bar \ep$ such that
$K+B$ has a bounded $n$-complement with $I\mid n$
if there are two divisors $E$
\textup(exceptional or not\textup) with discrepancy 
$a(E,X,B)\le -1+\bar \ep$\textup; this
$\bar \ep$ depends only on $\R$.
\end{corollary}

\begin{corollary}[cf. {\cite{Blache-1995}}, {\cite{Shokurov-2000}}]
\label{main-result-3-Cor-1-aC}
Fix a finite rational subset $\R\subset [0,\, 1]\cap \QQ$.
Let $(X,D)$ be a two-dimensional $0$-pair such that
the multiplicities of $B$ are contained in $\Phi(\R)$.
Then there exists a positive integer $n$ such that
$n(K+D)\sim 0$\textup; this
$n$ depends only on $\R$.
\end{corollary}

We give a sketch of the proof of our main results in 
Section \ref{sect-reduction}. One can see that 
our proof essentially uses reduction to lower-dimensional 
global pairs.
However it is expected that an improvement of our method
can use reduction to local questions in the same dimension. In fact 
we hope that the hypothesis in 
our main theorem \ref{main-result} should be 
the existence of local complements and 
Conjecture~\xref{BAB} for $\ep$-lt Fano varieties 
(without a boundary), where $\ep\ge \ep^o>0$,
$\ep^o$ is a constant depending only on the dimension
(cf. Addendum \ref{addendum-1}
and Corollary \ref{main-result-3-Cor-1-b}).
If $\dim X=2$, we can take $\ep^o=1/7$.
Note also that our main theorem \ref{main-result}
is weaker than one can expect. 
We think that the pair $(X,B)$ 
can be taken arbitrary log-semi-Fano (possibly
not klt and not FT) and possible boundary multiplicities 
can be taken arbitrary
real numbers in $[0,1]$ (not only in $\Phi(\R)$). 
The only hypothesis we have to assume is the existence of 
an $\RR$-complement $B^+\ge B$ (cf. \cite{Shokurov-2000}).
However the general case needs actually a {\em finite}
set of natural numbers for complements, and
there are no such universal number for all 
complements (cf. \cite[Example 5.2.1]{Shokurov-1992-e}). 

\subsection*{Acknowledgements}
The work was conceived 
in 2000 when the first 
author visited the Johns Hopkins University and 
finished during his stay 
in Max-Planck-Institut f\"ur Mathematik, Bonn
in 2006. 
He would like to thank these institutes 
for hospitality.
Finally both authors are
grateful to the referee
whose constructive criticism 
helped us to revise the paper very much.

\section{Preliminaries}
\subsection{Notation}
All varieties are assumed to be algebraic and defined over
an algebraically closed field $\Bbbk$ of characteristic zero.
Actually, main results
holds for any $\Bbbk$ of characteristic zero not
necessarily algebraically closed
since they are related to singularities
of general members of linear systems
(see \cite[5.1]{Shokurov-1992-e}). 
We use standard terminology and notation of the
Log Minimal Model Program (LMMP) \cite{KMM}, \cite{Utah},
\cite{Shokurov-1992-e}. For the definition of complements and their properties
we refer to \cite{Shokurov-1992-e}, \cite{Shokurov-2000},
\cite{Prokhorov-2001} and \cite{Prokhorov-Shokurov-2001}.
Recall that a \emph{log pair} (or a \emph{log variety}) is a pair
$(X,D)$ consisting of a normal variety $X$ and a \emph{boundary}
$D$, i.e., an $\RR$-divisor $D=\sum d_i D_i$ with multiplicities $0\le
d_i\le 1$.
As usual $K_X$ denotes the canonical (Weil) divisor of 
a variety $X$. Sometimes we will write $K$ instead of $K_X$ if
no confusion is likely. 
Everywhere below $a(E,X,D)$ denotes the discrepancy of $E$ with
respect to $K_X+D$. 
Recall the standard notation:
\[
\begin{array}{lll}
\discr(X,D)&=&\inf_E\{ a(E,X,D)\mid \codim\Center_X(E)\ge 2\},
\vspace{6pt}\\
\totaldiscr(X,D)&=&\inf_E\{ a(E,X,D)\mid
\codim\Center_X(E)\ge 1\}.
\end{array}
\]
In the paper we use the following strong
version of $\ep$-log terminal and $\ep$-log canonical
property. 
\begin{definition}
\label{def-ep-lt}
A log pair $(X,B)$ is said to be
\textit{$\ep$-log terminal} (\textit{$\ep$-log canonical})
if $\totaldiscr(X,B)>-1+\ep$ (resp., $\totaldiscr(X,B)\ge-1+\ep$).
\end{definition}

\subsection{}
Usually we work with $\RR$-divisors.
An $\RR$-divisor is an $\RR$-linear combination of 
prime Weil divisors.
An $\RR$-linear combination $D=\sum \alpha_i L_i$,
where the $L_i$ are integral Cartier divisors 
is called an $\RR$-Cartier divisor.
The pull-back $f^*$ of an $\RR$-Cartier divisor 
$D=\sum \alpha_i L_i$ under a morphism $f\colon Y\to X$
is defined as $f^*D:=\sum \alpha_i f^*L_i$.
Two $\RR$-divisors $D$ and $D'$ are said to be 
\emph{$\QQ$- \textup(resp., $\RR$-\textup)linearly equivalent}
if $D-D'$ is a
$\QQ$- (resp., $\RR$-)linear combination of principal divisors.
For a positive integer $I$, two $\RR$-divisors $D$ and $D'$ are said to be
\emph{$I$-linearly equivalent} if
$I(D-D')$ is an (integral) principal divisor.
The $\QQ$-linear (resp., $\RR$-linear, $I$-linear)
equivalence is denoted by $\q{\QQ}$
(resp., $\q{\RR}$, $\q{I}$).
Let $\Phi\subset \RR$ and let $D=\sum d_iD_i$ be an
$\RR$-divisor. We say that $D\in \Phi$ if $d_i\in \Phi$ for all
$i$.

\subsection{}
Let $f\colon X\to Z$ be a morphism of normal varieties.
For any $\RR$-divisor $\Delta$ on $Z$, define its 
\textit{divisorial pull-back} 
$f^\bullet\Delta$ as the closure of the usual pull-back $f^*\Delta$
over $Z\setminus V$, where $V$ is a closed subset 
of codimension $\ge 2$ such that $V\supset \Sing Z$
and $f$ is equidimensional over $Z\setminus V$. 
Thus each component of $f^\bullet\Delta$ dominates 
a component of $\Delta$.
It is easy to see that 
the divisorial pull-back $f^\bullet\Delta$ does not depend
on the choice of $V$.
Note however that in general $f^\bullet$ 
does not coincide with the usual
pull-back $f^*$ of $\RR$-Cartier divisors.

\begin{definition}
\label{definition}
Let $(X,B)$ be a log pair of global type (the latter means that
$X$ is projective). Then it is said to be
\begin{itemize}
\item[]
\emph{log Fano} variety if $K+B$ is lc and $-(K+B)$ is ample;
\item[]
\emph{weak log Fano} (WLF) variety if $K+B$ is lc and $-(K+B)$ is nef
and big;
\item[]
\emph{log semi-Fano} (ls-Fano) variety if $K+B$ is lc and $-(K+B)$ is nef;
\item[]
\emph{$0$-log pair} if $K+B$ is lc and numerically
trivial\footnote{Such a log pair can be called also a \emph{log
Calabi-Yau variety}. However the last notion usually assumes some
additional conditions such as $\pi_1(X)=0$ or $q(X)=0$.}.
\end{itemize}
In dimension two we usually use the word \textit{del Pezzo}
instead of \textit{Fano}.
\end{definition}

\begin{lemma-definition}
\label{CY->F}
Let $X$ be a normal projective variety. We say that
$X$ is \emph{FT} \textup(\emph{Fano type}\textup) if it
satisfies the following equivalent
conditions:
\begin{enumerate}
\item
there is a $\QQ$-boundary $\Xi$ such that $(X,\Xi)$ is a klt log Fano;
\item
there is a $\QQ$-boundary $\Xi$ such that $(X,\Xi)$ is a klt weak log Fano;
\item
there is a $\QQ$-boundary $\Theta$ such that $(X,\Theta)$ is a klt $0$-pair
and the components of $\Theta$ generate $N^1(X)$;
\item
for any divisor $\Upsilon$ there is
a $\QQ$-boundary $\Theta$ such that $(X,\Theta)$ is a klt $0$-pair and
$\Supp \Upsilon\subset \Supp \Theta$.
\end{enumerate}
\end{lemma-definition}
Similarly one can define relative FT and $0$-varieties $X/Z$,
and the results below hold for them too.

\begin{proof}
Implications (i) $\Longrightarrow$ (iv), (iv) $\Longrightarrow$ (iii), (i) $\Longrightarrow$ (ii)
are obvious and (ii) $\Longrightarrow$ (i)
follows by Kodaira's lemma (see, e.g., \cite[Lemma 0-3-3]{KMM}).
We prove (iii) $\Longrightarrow$ (i).
Let $(X,\Theta)$ be such as in (iii).
Take an ample divisor $H$ such that
$\Supp{H}\subset\Supp{\Theta}$ and put
$\Xi=\Theta-\ep H$, for $0<\ep\ll 1$.
Clearly, $(X,\Xi)$ is a klt log Fano.
\end{proof}

Recall that for any \textup(not necessarily effective\textup)
$\RR$-divisor $D$ on a variety $X$ a $D$-MMP is a sequence 
$X=X_1 \dashrightarrow X_N$
of extremal $D$-negative divisorial contractions and $D$-flips 
which terminates on a variety $X_N$ where either 
the proper transform of $D$ is 
nef or there exists
a $D$-negative contraction to a lower-dimensional variety 
\textup(see \cite[2.26]{Utah}\textup).

\begin{corollary}
\label{CY-MMP}
Let $X$ be an FT variety. Assume the LMMP in dimension $\dim X$.
Then the $D$-MMP works on $X$ with respect to any
$\RR$-divisor $D$.
\end{corollary}
\begin{proof}
Immediate by Lemma \ref{CY->F}, (iv).
Indeed, in the above notation we may assume that 
$\Supp D\subset\Supp \Theta$.
It remains to note that 
the $D$-LMMP is is nothing but the LMMP with respect to $(X,\Theta+\ep D)$
some $0<\ep\ll 1$. 
\end{proof}

\begin{lemma}
\label{lemma-FT}
\begin{enumerate}
\item
Let $f\colon X\to Z$ be a \textup(not necessarily birational\textup) 
contraction of normal varieties.
If $X$ is FT, then so is $Z$. 
\item
The FT property is preserved under birational
divisorial contractions and flips.
\item
Let $(X,D)$ be an ls-Fano variety such that $X$ is FT.
Let $f\colon Y\to X$ be a birational extraction such that
$a(E,X,D)< 0$ for every $f$-exceptional divisor $E$ over $X$.
Then $Y$ is also FT.
\end{enumerate}
\end{lemma}

We need the following result of Ambro
\cite[Th. 0.2]{Ambro-2005}
which is a variant of Log Canonical Adjunction
(cf. \ref{conj-main-adj}, \cite{Fujino-1999app}).
\begin{theorem}[{\cite[Th. 0.2]{Ambro-2005}}]
\label{th-ambro}
Let $(X,D)$ be a projective klt log pair, let
$f\colon X\to Z$ be a contraction, 
and let $L$ be a $\QQ$-Cartier 
divisor on $Z$ such that
\[
K+D\q{\QQ} f^*L.
\]
Then there exists a $\QQ$-Weil divisor $D_Z$ such 
that $(Z,D_Z)$ is a log variety with Kawamata log
terminal singularities and 
$L \q{\QQ} K_Z+D_Z$.
\end{theorem}

\begin{proof}[Proof of Lemma \xref{lemma-FT}]
First note that 
(ii) and the birational case of (i) easily follows from 
from \ref{CY->F} (iii). 
To prove (i) in the general case 
we apply Theorem \ref{th-ambro}.
Let $\Theta=\sum_i \theta_i\Theta_i$
be a $\QQ$-boundary on $X$ whose components generate $N^1(X)$
and such that $(X, \Theta)$ is a klt $0$-pair.
Let $A$ be an ample divisor on $Z$. By our assumption
$f^*A\equiv \sum_i \delta_i \Theta_i$.
Take $0<\delta\ll 1$ and put
$\Theta':= \sum_i (\theta_i-\delta \delta_i) \Theta_i$.
Clearly, $K+\Theta'\equiv -\delta f^*A$
and $(X,\Theta')$ is a klt log semi-Fano variety.
By the base point free theorem 
$K+\Theta'\q{\QQ} -\delta f^*A$.
Now by Theorem \ref{th-ambro} there is a $\QQ$-boundary
$\Theta_Z$ such that $(Z,\Theta_Z)$ is klt and 
$K_Z+\Theta_Z \q{\QQ} -\delta A$.
Hence $(Z, \Theta_Z)$ is a klt log Fano variety.
This proves (i).

Now we prove (iii). Let $\Xi$ be a boundary
such that $(X,\Xi)$ is a klt log Fano.
Let $D_Y$ and $\Xi_Y$ be proper transforms of $D$ and $\Xi$,
respectively. Then $(Y,D_Y)$ is an
ls-Fano, $(Y,\Xi_Y)$ is klt and $-(K_Y+\Xi_Y)$ is nef and big.
However $\Xi_Y$ is not necessarily a boundary.
To improve the situation we put
$\Xi':=(1-\ep)D_Y+\ep \Xi_Y$ for small positive $\ep$.
Then $(Y,\Xi')$ is a klt weak log Fano.
\end{proof}

\begin{definition}
\label{def-compl}
Let $X$ be a normal variety and
let $D$ be an $\RR$-divisor on $X$. Then a
\textit{$\QQ$-complement} of $K_X+D$ is a
log divisor $K_X+D'$ such that $D'\ge D$, $K_X+D'$ is lc and
$n(K_X+D')\sim 0$ for some positive integer $n$.

Now let $D=S+B$, where 
$B$ and $S$ have no common components, $S$ is an effective
integral divisor and $\down{B}\le 0$. Then we say that $K_X+D$ is
\textit{$n$-complemented}, if there is a $\QQ$-divisor $D^+$ such
that
\begin{enumerate}
\item
$n(K_X+D^+)\sim 0$ (in particular, $nD^+$ is integral divisor);
\item
$K_X+D^+$ is lc;
\item
$nD^+\ge nS+\down{(n+1)B}$.
\end{enumerate}
In this situation, $K_X+D^+$ is called an
\textit{$n$-complement} of $K_X+D$.
\end{definition}
Note that an $n$-complement is not necessarily 
a $\QQ$-complement (cf. Lemma \ref{lemma-PPP-n-1}).
\begin{remark}
\label{rem-complne}
Under (i) and (ii) of \ref{def-compl}, the condition (iii) follows from the 
inequality $D^+\ge D$. Indeed, write $D=\sum d_iD_i$ and $D^+=\sum d_i^+D_i$.
We may assume that $d_i^+<1$. Then we have 
\[
 nd_i^+\ge \down{nd_i^+}=\down{(n+1)d_i^+}\ge \down{(n+1)d_i}.
\]
\end{remark}

\begin{corollary}
\label{cor-complne}
Let $D^+$ be an $n$-complement of $D$ such that $D^+\ge D$.
Then $D^+$ is also an $nI$-complement of $D$ for any
positive integer $I$.
\end{corollary}

For basic properties of complements we refer to 
\cite[\S 5]{Shokurov-1992-e} and \cite{Prokhorov-2001}, see also 
\S \ref{sect-hyper}.

\subsection{}
Fix a class of (relative) log pairs $(\mathcal{X}/\mathcal{Z}\ni o,\,
\mathcal{B})$, where $o$ is a point on each $Z\in \mathcal Z$.
We say that this class has \emph{bounded complements} if
there is a constant $\operatorname{Const}$ such that for any log
pair $(X/Z,B)\in (\mathcal{X}/\mathcal{Z},\,
\mathcal{B})$ the log divisor $K+B$ is
$n$-complemented near the fibre over $o$ for some
$n\le\operatorname{Const}$.

\subsection{Notation}
Let $X$ be a normal $d$-dimensional variety and let
$\BBB=\sum_{i=1}^r B_i$ be any reduced divisor on $X$. Recall that
$Z_{d-1}(X)$ usually denotes the group of Weil divisors on $X$.
Consider the vector space $\mathfrak{D}_{\BBB}$ of all
$\RR$-divisors supported in $\BBB$:
\[
\mathfrak{D}_{\BBB}:=\bigl\{D\in Z_{d-1}(X)\otimes \RR\mid
\Supp{D}\subset\BBB\bigr\}=\sum_{i=1}^r \RR\cdot B_i.
\]
As usual, define a norm in $\mathfrak{D}_{\BBB}$ by
\[
\|B \|=\max(|b_1|,\dots, |b_r|),
\]
where $B=\sum_{i=1}^rb_iB_i\in \mathfrak{D}_{\BBB}$. For any
$\RR$-divisor $B=\sum_{i=1}^r b_i B_i$, put
$\mathfrak{D}_{B}:=\mathfrak{D}_{\Supp{B}}$.

\section{Hyperstandard multiplicities}
\label{sect-hyper}
Recall that \textit{standard multiplicities} $1-1/m$ naturally
appear as multiplicities in the divisorial adjunction formula
$(K_X+S)|_S=K_S+\Diff_S$ (see \cite[\S 3]{Shokurov-1992-e},
\cite[Ch. 16]{Utah}). Considering the adjunction formula for fibre
spaces and adjunction for higher codimensional 
subvarieties one needs to introduce a
bigger class of multiplicities.

\begin{example}
\label{ex-Kodaira}
Let $f\colon X\to Z\ni P$ be a minimal two-dimensional elliptic
fibration over a one-dimensional germ
($X$ is smooth). We can write a natural formula
$K_X=f^*(K_Z+D\sdiv)$, where $D\sdiv=d_P P$ is an effective divisor
(cf. \ref{def-cW} below).
From Kodaira's classification of singular fibres
(see \cite{Kodaira-1963}) we obtain the following
values of $d_P$:
\par\medskip\noindent
\begin{center}
\begin{tabular}{p{35pt}|p{35pt}p{25pt}p{25pt}
p{25pt}p{25pt}p{25pt}p{25pt}p{25pt}} {\rm Type}&
$m\mathrm{I}_n$&$\mathrm{II}$&$\mathrm{III}$&$\mathrm{IV}$&
$\mathrm{I}^*_b$&$\mathrm{II}^*$&$\mathrm{III}^*$&$\mathrm{IV}^*$
\vspace{10pt}
\\
\hline
\\
$d_P$& $1-\frac1m$&$\frac16$&$\frac14$&$\frac13$&$\frac12$
&$\frac56$&$\frac34$&$\frac23$\\
\end{tabular}
\end{center}
\medskip
Thus the multiplicities of $D\sdiv$ are not necessarily standard.
\end{example}

\subsection{}
\label{not-hyperstand}
Fix a subset $\R\subset \RR_{\ge 0}$.
Define
\[
\Phi(\R):=\left\{\left. 1-\frac rm \quad \right|\quad m\in \ZZ,
\quad m>0 \quad
r\in \R \right\} \bigcap \Bigl[0,\, 1\Bigr].
\]
We say that an $\RR$-boundary $B$ has \textit{hyperstandard multiplicities}
with respect to $\R$ if $B\in \Phi(\R)$. For example, if
$\R=\{0,1\}$, then $\Phi(\R)$ is the set of \textit{standard multiplicities}.
The set $\R$ is said to be \emph{rational} if $\R\subset \QQ$.
Usually we will assume that $\R$ is rational and finite.
In this case we denote
\[
I(\R):=
\lcm\bigl(
\text{denominators of $r\in\R\setminus \{0\}$}
\bigr).
\]

\para
\label{para-def-N}
Denote by
$\NNN_d(\R)$ the set
of all $m\in \ZZ$, $m>0$ such that there exists
a log semi-Fano variety $(X,D)$ of dimension $\le d$ satisfying
the following properties:
\begin{enumerate}
\item
$X$ is FT and $D\in\Phi(\R)$;
\item
either $(X,D)$ is klt or $K_X+D\equiv 0$;
\item
$K_X+D$ is $m$-complemented, $I(\R) \mid m$, and $m$ is
minimal under these conditions.
\end{enumerate}
Since any nef divisor is semiample on FT
variety, for any log semi-Fano variety $(X,D)$ 
satisfying (i) and (ii),
there exists some $m$ in (iii). 
Put
\[
N_{d}=N_{d}(\R):=\sup \NNN_d(\R), \quad
\ep_{d}=\ep_{d}(\R):=1/(N_{d}+2).
\]
We expect that $\NNN_d(\R)$ is bounded whenever
$\R$ is finite and rational, see
Theorems \ref{main-result} and \ref{main-result0}. In particular,
$N_{d}<\infty$ and $\ep_{d}>0$.
For $\ep\ge 0$, define also the set of
\textit{semi-hyperstandard} multiplicities
\[
\Phi(\R,\ep):=
\Phi(\R) \cup [1-\ep,\, 1].
\]

Fix a positive integer $n$ and define the set $\PPP_n\subset \RR$ by
\[
\label{def_PPP_n}
\alpha\in\PPP_n\quad\Longleftrightarrow\quad 0\le\alpha\le
1\quad\text{and}\quad \down{(n+1)\alpha}\ge n\alpha.
\index{$\PPP_n$}
\]
This set obviously satisfies the following property:

\begin{lemma}
\label{lemma-PPP-n-1}
If $D\in \PPP_n$ and $D^+$ is an $n$-complement, then
$D^+\ge D$.
\end{lemma}

Taking \ref{cor-complne} into account we immediately obtain
the following important.
\begin{corollary}
\label{cor-imp}
Let $D\in \PPP_n$ and let $D^+$ be an $n$-complement 
of $D$. Then $D^+$ is an $nI$-complement of $D$ for any 
positive integer $I$.
\end{corollary}

\begin{lemma}[cf. {\cite[Lemma 2.7]{Shokurov-2000}}]
\label{lemma-PPP-n}
If $\R\subset [0,\, 1]\cap \QQ$, $I(\R) \mid n$,
and $0\le \ep \le 1/(n+1)$, then
\[
\PPP_n\supset \Phi (\R,\ep).
\]
\end{lemma}

\begin{proof}
Let $1\ge \alpha\in \Phi (\R,\ep)$. If $\alpha\ge 1-\ep$,
then
\[
(n+1)\alpha > n+1-\ep(n+1)\ge n.
\]
Hence, $\down{(n+1)\alpha}\ge n\ge n\alpha$ and $\alpha\in \PPP_n$.
Thus we may assume that $\alpha\in \Phi(\R)$.
It is sufficient to show that
\begin{equation}
\label{(n+1)}
\down{(n+1)\left(1-\frac rm\right)}\ge n\left(1-\frac rm\right)
\end{equation}
for all $r\in\R$ and $m\in\ZZ$, $m>0$. We may assume that $r>0$. It is
clear that \eqref{(n+1)} is equivalent to the following inequality
\begin{equation}
\label{(n+1)1}
(n+1)\left(1-\frac rm\right)\ge k\ge n\left(1-\frac rm\right),
\end{equation}
for some $k\in\ZZ$ (in fact, $k=\down{(n+1)\left(1-\frac
rm\right)}$). By our conditions, $N:=nr\in\ZZ$, $N>0$. Thus
\eqref{(n+1)1} can be rewritten as follows
\begin{equation}
\label{(n+1)2}
mn-N+m-r\ge mk\ge mn-N.
\end{equation}
Since $m-r\ge m-1$, inequality \eqref{(n+1)2} has a solution in
$k\in \ZZ$. This proves the statement.
\end{proof}

\begin{proposition}[{\cite[Prop. 4.3.2]{Prokhorov-2001}}, {\cite[Prop. 6.1]{Prokhorov-Shokurov-2001}}]
\label{cor-pull-back_compl-I}
Let $f\colon Y\to X$ be a birational contraction 
and let $D$ be an $\RR$-divisor on $Y$ such that
\begin{enumerate}
\item
$K_Y+D$ is nef over $X$,
\item
$f_*D\in \PPP_n$ \textup(in particular, $f_*D$ is a boundary\textup).
\end{enumerate}
Assume that $K_X+f_*D$ is $n$-complemented.
Then so is $K_Y+D$.
\end{proposition}

\begin{proposition}[{\cite[Prop. 4.4.1]{Prokhorov-2001}}, {\cite[Prop. 6.2]{Prokhorov-Shokurov-2001}}]
\label{prodolj}
Let $(X/Z\ni o,D=S+B)$ be a log variety. Set
$S:=\down{D}$ and $B:=\fr{D}$. Assume that
\begin{enumerate}
\item
$K_X+D$ is plt\textup;
\item
$-(K_X+D)$ is nef and big over $Z$\textup;
\item
$S\ne 0$ near $f^{-1}(o)$\textup;
\item
\label{prop-2-iv}
$D\in \PPP_n$.
\end{enumerate}

Further, assume that near $f^{-1}(o)\cap S$ there exists an
$n$-complement $K_S+\Diff_ S(B)^+$ of $K_S+\Diff_ S(B)$. Then near
$f^{-1}(o)$ there exists an $n$-complement $K_X+S+B^+$ of
$K_X+S+B$ such that $\Diff_S(B)^+=\Diff_S({B^+})$.
\end{proposition}

\subsection*{Adjunction on divisors (cf. \cite[Cor. 3.10, Lemma 4.2]{Shokurov-1992-e})}
Fix a subset $\R\subset \RR_{\ge 0}$. Define also the new set
\[
\ov \R:= \left\{\left. r_0-m \sum_{i=1}^s (1-r_i)\ \right|\
r_0,\dots, r_s\in\R, \ m\in \ZZ, \ m>0 \right\}\cap \RR_{\ge 0}.
\]
It is easy to see that $\ov \R\supset \R$. For example, if
$\R=\{0,1\}$, then $\ov\R=\R$.

\begin{lemma}
\begin{enumerate}
\item
If $\R\subset [0,\, 1]$, then $\ov \R\subset [0,\, 1]$. 
\item
If $\R$ is finite and rational, then so is $\ov \R$.
\item
$I(\R)=I(\ov \R)$.
\item
Let $\mathfrak{G}\subset \QQ$ be an additive subgroup containing
$1$ and let $\R=\R_{\mathfrak{G}} :=\mathfrak{G}\cap [0,1]$. Then
$\ov \R=\R$.
\item
If the ascending chain condition \textup(a.c.c.\textup) holds for the set $\R$,
then it holds for $\ov \R$.
\end{enumerate}
\end{lemma}

\begin{proof}
(i)-(iv) are obvious. We prove (v). Indeed,
let
\[
q^{(n)}=r_0^{(n)}-m^{(n)} \sum_{i=1}^{s^{(n)}} (1-r_i^{(n)})\in
\ov\R
\]
be an infinite increasing sequence, where $r_i^{(n)}\in\R$ and
$m^{(n)}\in\ZZ_{>0}$. By passing to a subsequence, we may assume that
$m^{(n)} \sum_{i=1}^{s^{(n)}} (1-r_i^{(n)})>0$, in particular, $s^{(n)}>0$
for all $n$. There is a constant $\ep=\ep(\R)>0$ such that
$1-r_i^{(n)}>\ep$ whenever $r_i^{(n)}\neq 1$. Thus, $0\le
q^{(n)}\le r_0^{(n)}- m^{(n)}s^{(n)}\ep$ and $m^{(n)}s^{(n)}\le
(r_0^{(n)}- q^{(n)})/\ep$. Again by passing to a subsequence, we
may assume that $m^{(n)}$ and $s^{(n)}$ are constants:
$m^{(n)}=m$, $s^{(n)}=s$. Since the numbers $r_i^{(n)}$ satisfy a.c.c., 
the sequence
\[
q^{(n)}=r_0^{(n)}+m \sum_{i=1}^{s} r_i^{(n)}-ms
\]
is not increasing, a contradiction.
\end{proof}

\begin{proposition}
\label{prop-prop}
Let $\R\subset [0,\, 1]$, $1\in \R$, $\ep \in [0,\, 1]$, 
and let $(X,S+B)$ be a plt log
pair, where $S$ is a prime divisor, $B\ge 0$, and $\down B=0$. If
$B\in \Phi(\R, \ep)$, then $\Diff_S(B)\in \Phi(\ov \R, \ep)$.
\end{proposition}
\begin{proof}
Write $B=\sum b_i B_i$,
where the $B_i$ are prime divisors and $b_i \in \Phi(\R, \ep)$.
Let $V\subset S$ be a prime divisor. By \cite[Cor. 3.10]{Shokurov-1992-e}
the multiplicity $d$ of $\Diff_S(B)$ along $V$ is computed using the
following relation:
\begin{equation*}
\label{eq-adji-div-coeff}
d=1-\frac1{n}+\frac1n\sum_{i=0}^s k_i b_i=1-\frac{\beta}{n},
\end{equation*}
where $n$, $k_i\in \ZZ_{\ge 0}$, and $\beta:=1-\sum k_i b_i$.
It is easy to see that $d\ge b_i$ whenever $k_i>0$.
If $b_i\ge 1-\ep$, this implies $d\ge 1-\ep$.
Thus we may assume that 
$b_i\in \Phi(\R)$ whenever $k_i>0$. Therefore,
\[
\beta=1-\sum k_i\left(1-\frac{r_i}{m_i}\right),
\]
where $m_i\in\ZZ_{>0}$, $r_i\in \R$. Since $(X,S+B)$ is plt, $d<1$.
Hence, $\beta>0$. 
If $m_i=1$ for all $i$, then 
\[
\beta=1- \sum k_i(1-r_i)\in \ov \R.
\]
So, $d\in \Phi(\ov \R)$
in this case. Thus we may assume that $m_0>1$. 
Since
$1-\frac{r_i}{m_i}\ge 1-\frac{1}{m_i}$, we have $m_1=\cdots=m_s=1$ and $k_0=1$. 
Thus,
\[
\beta= \frac{r_0}{m_0}-
\sum _{i=1}^s k_i(1-r_i)= \frac{r_0-m_0\sum_{i=1}^sk_i(1-r_i)}{m_0}
\]
and $m_0\beta=r_0-m_0\sum_{i=1}^sk_i(1-r_i)\in \ov \R$.
Hence, $d=1-\frac{m_0\beta}{m_0n}\in \Phi(\ov \R)$.
\end{proof}

\begin{proposition}
\label{lemma-compl-P1-bound-s}
Let $1\in \R\subset [0,\, 1]$ and 
let $(X,B)$ be a klt log semi-Fano of dimension $\le d$
such that $X$ is FT.
Assume the LMMP in dimension $d$.
If $B\in \Phi(\R,\ep_{d})$, then
there is an $n$-complement $K+B^+$ of $K+B$ for some $n\in \NNN_d(\R)$.
Moreover, $B\in \PPP_n$, and so $B^+\ge B$. 
\end{proposition}
In the proposition we do not assert that
$\NNN_d(\R)$ is finite.
However later on we use the proposition in the induction process
when the set of indices is finite (cf. the proof
of Lemma \ref{lemma-2-compl-} and see \ref{subsectio-reduction}
-- \ref{sufficient-condition}). 
\begin{proof}
If $\ep _{d}=0$, then $\Phi(\R,\ep_{d})=\Phi(\R)$
and there is nothing to prove.
So we assume that $\ep_{d}>0$.
If $X$ is not $\QQ$-factorial, we
replace $X$ with its small $\QQ$-factorial modification.
Write $B=\sum b_iB_i$. Consider the new boundary
$D=\sum d_i B_i$, where
\[
d_i=
\begin{cases}
b_i&\text{if $b_i< 1-\ep _{d}$,}
\\
1-\ep _{d}&\text{otherwise.}
\end{cases}
\]
Clearly, $D\in \Phi(\R)$. Since $D\le B$, there is
a klt $\QQ$-complement $K+D+\Lambda$ of 
$K+D$ (by definition, $\Lambda\ge 0$).
Run $-(K+D)$-MMP. Since all the birational transformations 
are $K+D+\Lambda$-crepant, 
they preserve the klt property of
$(X,D+\Lambda)$ and $(X,D)$. 
Each extremal ray is $\Lambda$-negative, and therefore is birational.
At the end we get a model $(\bar X,\bar D)$ which is a log semi-Fano
variety.
By definition, since $\bar D\in \Phi(\R)$ and $X$ is FT, there is an
$n$-complement $\bar D^+$ of $K_{\bar X}+\bar D$ for some $n\in \NNN_d(\R)$.
Note that 
$I(\R) \mid n$, so 
$\bar D\in \Phi(\R)\subset \PPP_n$.
By Proposition \ref{cor-pull-back_compl-I} we can pull-back this complement
to $X$ and this gives us an $n$-complement of $K_X+B$.
The last assertion follows by Lemmas \ref{lemma-PPP-n} and
\ref{lemma-PPP-n-1}.
\end{proof}

\section{General reduction}
\label{sect-reduction}
In this section we outline the 
main reduction step in the proof of our main results
\ref{main-result} and \ref{main-result0}.
First we concentrate on the klt case.
The non-klt case of \ref{main-result0}
will be treated in 
\ref{subsect-1} and \ref{subsetion-0-pairs-proof}.
Note that in Theorem
\ref{main-result} it is sufficient to find only 
one integer $n=n(d,\R)$ divisible by $I(\R)$ and such that 
$K+B$ is $n$-complemented. Other statements 
immediately follows by 
Corollary \ref{cor-imp},
Lemma \ref{lemma-PPP-n} and by the definition of complements.

\subsection{Setup}
\label{eq-reduction-setup}
Let $(X,B=\sum b_iB_i)$ be a klt log semi-Fano variety of dimension $d$
such that $B\in \Phi(\R)$.
In particular, $B$ is a $\QQ$-divisor.
Assume that $X$ is FT. 
By induction we may assume 
that Theorems 
\ref{main-result} and \ref{main-result0} hold in dimension $d-1$. 
So, by this inductive hypothesis, $\ep_{d-1}(\ov\R)>0$ whenever
$\R\subset [0,\, 1]$ is finite and rational. Take any $0<\ep'\le \ep_{d-1}(\ov\R)$.
Put also $I:=I(\R)$.

\para{}
\label{para-rem-Noetherian}
First assume that the pair $(X,B)$ is
$\ep'$-lt. Then the multiplicities of
$B$ are contained in the finite set $\Phi(\R)\cap [0,\, 1-\ep']$.
By Conjecture \xref{BAB} the pair $(X,B)$ is bounded.
Hence $(X,\Supp B)$ belongs to an algebraic family and we may assume that 
the multiplicities of $B$ are fixed. Let $m:=nI$.
The condition that $K+B$ is $m$-complemented is equivalent to the following 
\[
\exists\ \ov B\in \left|-K-\down{(m+1)B}\right|\hspace{4pt} \text{such that}
\hspace{4pt}
\left(X, \frac 1m \left(\down{(m+1)B}+\bar B\right)\right)\hspace{4pt} \text{is lc}
\]
(see \ref{def-compl}, \cite[5.1]{Shokurov-1992-e}).
Obviously, the last condition is open in the deformation space of 
$(X,\Supp B)$.
By Proposition \ref{prop-1} below and Noetherian induction 
the log divisor $K+B$ has a bounded $nI$-complement
for some $n\le C(d,\R)$.
From now on we assume that $(X,B)$ is not
$\ep'$-lt.

\subsection{}
\label{sub-reduction-1}
We replace $(X,B)$ with log crepant $\QQ$-factorial
blowup of all divisors $E$ of discrepancy
$a(E,X,B)\le -1+\ep'$, see \cite[21.6.1]{Utah}.
Condition $B\in \Phi(\R)$ will be replaced with $B\in \Phi(\R,\ep')\cap \QQ$.
Note that our new $X$ is again FT by Lemma \ref{lemma-FT}.
From now on we assume that $X$ is $\QQ$-factorial and
\begin{equation}
\label{eq-discr-2-B-ep}
\discr(X,B)>-1+\ep'.
\end{equation}

\para{}
\label{reduct-FT-Lambda}
For some $n_0\gg 0$, the divisor
$n_0B$ is integral and the linear system $|-n_0(K+B)|$
is base point free. Let $\ov B\in |-n_0(K+B)|$ be a general member.
Put $\Theta:=B+\frac1 {n_0} \ov B$. By Bertini's theorem 
$\discr(X,\Theta)=\discr(X,B)$.
Thus we have the following
\begin{itemize}
\item 
$K_X+\Theta$ is a klt $\QQ$-complement of $K_X+B$,
\item
$\discr(X,\Theta)\ge 1+\ep'$, and
\item
$\Theta-B$ is supported in a movable (possibly trivial) divisor.
\end{itemize}
Define a new boundary $D$ with $\Supp D=\Supp B$:
\begin{equation}
\label{eq-D}
D:=\sum d_iB_i,\quad\text{ where}\quad d_i=
\begin{cases}
1 & \text{if}\quad b_i\ge 1-\ep',
\\
b_i &
\text{otherwise}.
\end{cases}
\end{equation}
Here the $B_i$ are components of $B$.
Clearly, $D\in \Phi(\R)$, $D>B$, and by \eqref{eq-discr-2-B-ep}
we have $\down{D}\neq 0$.

\begin{lemma}[the simplest case of the global-to-local statement]
\label{lemma-2-compl-}
Fix a finite set $\R\subset [0,\, 1]\cap \QQ$.
Let $(X\ni o,D)$ be the germ of a $\QQ$-factorial klt
$d$-dimensional singularity, where $D\in \Phi (\R,\ep_{d-1}(\ov\R))$.
Then there is an $n$-complement of $K_X+D$ with
$n\in \NNN_{d-1}(\ov \R)$.
\end{lemma}
Recall that according to our inductive hypothesis,
$\NNN_{d-1}(\ov \R)$ is finite. 
\begin{proof}
Consider a plt blowup $f\colon \tilde X\to X$ of $(X,D)$
(see \cite[Prop. 3.6]{Prokhorov-Shokurov-2001}).
By definition the exceptional locus of $f$ is an
irreducible divisor $E$, the pair $(\tilde X, \tilde D+E)$ is plt, and
$-(K_{\tilde X}+\tilde D+E)$ is $f$-ample, where $\tilde D$
is the proper transform of $D$. We can take $f$ so that $f(E)=o$,
i.e., $E$ is projective. By Adjunction $-(K_E+ \Diff_E(\tilde D))$ is ample and
$(E, \Diff_E(\tilde D))$ is klt. By Proposition \ref{prop-prop} we have
$\Diff_E(\tilde D)\in \Phi(\ov \R, \ep_{d-1}(\ov\R))$.
Hence there is an $n$-complement of $K_E+\Diff_E(\tilde D)$ with
$n\in \NNN_{d-1}(\ov \R)$, see Proposition \ref{lemma-compl-P1-bound-s}.
This complement can be extended to $\tilde X$ by Proposition \ref{prodolj}.
\end{proof}

\begin{claim}
\label{claim-0}
The pair $(X,D)$ is lc.
\end{claim}
\begin{proof}
By Lemma \ref{lemma-2-compl-} near each point $P\in X$ there is
an $n$-complement $K+B^+$ of $K+B$ with $n\in \NNN_{d-1}(\ov \R)$.
By Lemma \ref{lemma-PPP-n}, we have $\PPP_n\supset \Phi (\R,\ep)$.
Hence, by Lemma \ref{lemma-PPP-n-1}, $B^+\ge B$. On the other hand,
$nB^+$ is integral and for any component of $D-B$, its multiplicity
in $B$ is $\ge \ep'>1/(n+1)$. 
Hence, $B^+\ge D$ and so $(X,D)$ is lc near $P$.
\end{proof}

\subsection{}
\label{subs-2-inductive-steps}
Run $-(K+D)$-MMP (anti-MMP).
If $X$ is FT, this is possible by Corollary \ref{CY-MMP}. 
Otherwise $K+B\equiv 0$ and $-(K+D)$-MMP coincides with
$K+B-\delta (D-B)$-MMP for some small positive $\delta$.

It is clear that
property $B\in \Phi(\R,\ep')$ is preserved on each step.
All birational transformations are $(K+\Theta)$-crepant.
Therefore $K+\Theta$ is klt on each step. Since $B\le \Theta$,
$K+B$ is also klt.
By Claim \ref{claim-0} the log canonical property of $(X,D)$ is also preserved
and $X$ is FT on each step by \ref{lemma-FT}.

\begin{claim1}
\label{red-not-contr}
None of components of $\down D$ is contracted.
\end{claim1}
\begin{proof}
Let $\var\colon X\to \bar X$ be a $K+D$-positive extremal
contraction and let $E$ be the corresponding exceptional divisor.
Assume that $E\subset \down D$. Put $\bar D:=\var_* D$.
Since $K_X+D$ is $\var$-ample, we can write
\[
K_X+D=\var^*(K_{\bar X}+\bar D)-\alpha E,\qquad \alpha> 0.
\]
On the other hand, since $(\bar X,\bar D)$ is lc, we have
\[
-1 \le a(E,\bar X,\bar D)=a(E,X,D)-\alpha=-1-\alpha<-1,
\]
a contradiction.
\end{proof}

\begin{corollary}
\label{cor-reduct-discrs}
Condition \eqref{eq-discr-2-B-ep} holds on each step of our MMP.
\end{corollary}
\begin{proof}
Note that all our birational transformations are $(K+\Theta)$-crepant.
Hence by \eqref{reduct-FT-Lambda}, it is sufficient to show that
none of the components of $\Theta$ of multiplicity $\ge 1-\ep'$
is contracted.
Assume that on some step we contract a component $B_i$
of multiplicity
$b_i\ge 1-\ep'$. Then by \eqref{eq-D} $B_i$ is a component of $\down D$.
This contradicts Claim \ref{red-not-contr}.
\end{proof}

\subsection{Reduction}
\label{subsectio-reduction}
After a number of divisorial contractions and flips
\begin{equation}
\label{MMP}
X\dashrightarrow X_1\dashrightarrow \cdots\dashrightarrow X_N=Y,
\end{equation}
we get a $\QQ$-factorial model $Y$ such that either
\para{}
\label{cases-1}
there is a non-birational
$K_{Y}+D_Y$-positive extremal contraction $\var\colon Y\to Z$ to a
lower-dimensional variety $Z$, or
\para{}
\label{cases-2}
$-(K_{Y}+D_Y)$ is nef.

Here $\square_Y$ denotes
the proper transform of $\square$ on $Y$.
\begin{claim}
\label{claim-1}
In case \eqref{cases-1}, $Z$ is a point, i.e., $\rho(Y)=1$
and $-(K_Y+B_Y)$ is nef.
\end{claim}
\begin{proof}
Let $F=\var^{-1}(o)$ be a
general fibre. Since $\rho(Y/Z)=1$ and $-(K+B)\equiv \Theta-B\ge 0$,
the restriction $-(K+B)|_F$ is nef.
It is clear that $B|_F\in \Phi(\R,\ep')$.
Assume that
$Z$ is of positive dimension. Then $\dim F<\dim X$.
By our inductive hypothesis and Proposition \ref{lemma-compl-P1-bound-s}
there is a bounded $n$-complement
$K_F+B|_F^+$ of $K_F+B|_F$ for some $n\in \NNN_{d-1}(\R)\subset \NNN_{d-1}(\ov\R)$.
By Lemmas \ref{lemma-PPP-n-1} and \ref{lemma-PPP-n},
we have $B|_F^+\ge D|_F \ge B|_F$. On the other hand,
$(K_X+D)|_F$ is $\var$-ample, a contradiction.
\end{proof}

\subsection{}
\label{case}
Therefore we have a $\QQ$-factorial FT variety $Y$ and two boundaries $B_Y=\sum b_iB_i$
and $D_Y=\sum d_i B_i$
such that $\discr(Y,B_Y)>-1+\ep'$, $B_Y\in \Phi(\R,\ep')$,
$D_Y\in \Phi(\R)$, $D_Y\ge B_Y$, and $d_i>b_i$ if and only if
$d_i=1$ and $b_i\ge -1+\ep'$.
Moreover, one of
the following two cases holds:
\para{}
\label{case-rho=1}
$\rho(Y)=1$, $K_Y+D_Y$ is ample, and $(Y,B_Y)$ is a klt log semi-Fano variety, or
\para{}
\label{case-nef}
$({Y},D_Y)$ is a log semi-Fano variety with $\down{D_Y}\neq 0$.
Since $D>B$, this case does not occur if $K+B\equiv 0$.

These two cases will be treated in sections \xref{sec-rho=1} and
\xref{sec-nef}, respectively.

\subsection{Outline of the proof of Theorem \ref{main-result}}
\label{sufficient-condition}
Now we sketch the basic idea in the proof of boundedness in case \eqref{case-rho=1}.
By \eqref{para-rem-Noetherian} we may assume 
that $(X,B)$ is not
$\ep'$-lt.
Apply constructions of \ref{sub-reduction-1},
\ref{subs-2-inductive-steps} and \ref{subsectio-reduction}.
Recall that on each step of \eqref{MMP} we contract an extremal ray which is
$(K+D)$-positive. By Proposition \xref{cor-pull-back_compl-I}
we can pull-back $n$-complements with
$n\in \NNN_d(\R)$ of $K_{Y}+D_Y$ to our original $X$. However it can happen
in case $\rho(Y)=1$ that $K_{Y}+D_Y$ has no
any complements. In this case we will show in
Section \ref{sec-rho=1} below that the multiplicities of $B_Y$ are
bounded from the above: $b_i<1-c$, where $c>0$. By Claim
\xref{red-not-contr} divisorial contractions in \eqref{MMP} do not
contracts components of $B$ with multiplicities $b_i\ge
1-\ep'$. Therefore the multiplicities of $B$ also are bounded
from the above. Combining this with $\discr (X,B)>-1+\ep'$
and Conjecture \xref{BAB} we get that $(X,\Supp B)$ belong to an
algebraic family. By Noetherian induction (cf. \eqref{para-rem-Noetherian}) 
we may assume that
$(X,\Supp B)$ is fixed. Finally, by Proposition \xref{prop-1} 
we have that $(X,B)$ has bounded
complements.

Case \eqref{case-nef} will be treated in Sect.~\xref{sec-nef}.
In fact in this case we study the contraction $f\colon Y\to Z$ 
given by $-(K+D)$. 
When $Z$ is a lower-dimensional variety,
$f$ is a fibration onto varieties 
with trivial log canonical divisor.
The existence of desired complements
can be established inductively, by using
an analog of Kodaira's canonical bundle formula
(see Conjecture \xref{conj-main-adj}).

The proof of Theorem \ref {main-result0} 
in case when $(X,B)$ is not klt is based on the following

\begin{lemma}
\label{lemma-Adj-hor-mult-1}
Let $(X,B=\sum b_iB_i)$ 
be a $0$-pair of dimension $d$ such that
$B\in \Phi(\R,\ep')$, where $\ep':=\ep_{d-1}(\ov \R)$. 
Assume the LMMP in dimension $d$.
Further, assume either 
\begin{enumerate}
\item 
$(X,B)$ is not klt and Theorems \xref{main-result}-\xref{main-result0}
hold in dimension $d-1$, or
\item 
Theorem \xref{main-result0}
holds in dimension $d$.
\end{enumerate}
Then there exists $\lambda:=\lambda(d,\R)>0$
such that either $b_i=1$ or $b_i\le 1-\lambda$
for all $b_i$. 
\end{lemma}
\begin{proof}
The proof is by induction on $d$.
Case $d=1$ is well-known;
see, e.g., Corollary \ref{lemma-compl-P1-bound}.
If $B=0$, there is nothing to prove.
So we assume that $B>0$.
Assume that, for some component 
$B_i=B_0$, we have $1-\lambda<b_0<1$.

First consider the case when $(X,B)$ is not klt.
Replace $(X,B)$ with its 
$\QQ$-factorial dlt modification so that $\down B\neq 0$.
If $B_0\cap \down B=\emptyset$, we run 
$K+B-b_0B_0$-MMP. After several steps 
we get a $0$-pair $(\hat X,\hat B)$ such that 
$\hat {\down B}\neq 0$, $\hat B_0\neq 0$ and one of the following holds:
\begin{enumerate}
\item[(a)]
$\hat B_0\cap \hat {\down B}\neq \emptyset$, or 
\item[(b)]
there is an extremal contraction $\hat f\colon \hat X\to \hat Z$
to a lower-dimensional variety 
such that both $\hat{\down B}$ and $\hat B_0$
are $\hat f$-ample, and $\hat B_0\cap \hat {\down B}= \emptyset$.
(In particular, $\hat Z$ is not a point.) 
\end{enumerate}
In the second case we can apply induction hypothesis 
restricting $\hat B$ to a general fibre of $\hat f$. In the first case,
replacing our original $(X,B)$ with a dlt modification
of $(\hat X,\hat B)$, we may assume that $B_0\cap \down B\neq \emptyset$.
Let $B_1\subset \down B$
be a component meeting $B_0$. 
Then $(B_1,\Diff_{B_1}(B-B_1))$ is a $0$-pair 
with $\Diff_{B_1}(B-B_1)\in \Phi(\ov \R,\ep')$
(see Lemma \ref{lemma-PPP-n}). 
Write $\Diff_{B_1}(B-B_1)=\sum \delta_i\Delta_i$
and let $\Delta_0$ be a component of $B_0\cap B_1$.
Then the multiplicity $\delta_0$ of $\Delta_0$ in $\Diff_{B_1}(B-B_1)$
is computed as follows: $\delta_0=1-1/r+\sum_l k_lb_l/r$, where 
$r$ and $k_l$ are non-negative integers and $r,\, k_0>0$ 
(see \cite[Corollary 3.10]{Shokurov-1992-e}). 
Since $b_0>1-\ep'>1/2$, we have 
$k_0=1$ and $k_l=0$ for $l\neq 0$.
Hence, $\delta_0=1-1/r+b_0/r$.
By induction
we may assume either $\delta_0=1$ or $\delta_0\le 1-\lambda(d-1,\ov R)$.
Thus we have either $b_0=1$ or 
\[
b_0\le 1-r\lambda(d-1,\ov R)\le 1-\lambda(d-1,\ov R)
\]
and we can put $\lambda(d, R)=\lambda(d-1,\ov R)$ in this case.

Now consider the case when $(X,B)$ is klt.
Replace $(X,B)$ with its 
$\QQ$-factorialization and again
run $K+B-b_0B_0$-MMP:
\ $(X,B) \dashrightarrow (X',B')$.
Clearly at the end we get a $B_0'$-positive extremal 
contraction $\varphi \colon (X',B')\to W$ to a lower-dimensional variety 
$W$. If $W$ is not a point, we can apply induction
restricting $B'$ to a general fibre. 
Thus replacing $(X,B)$ with $(X',B')$ we may assume that 
$X$ is $\QQ$-factorial, 
$\rho(X)=1$ and $(X,B-b_0 B_0)$ is a klt log Fano variety.
In particular, $X$ is FT.
By Theorem \ref{main-result0}
and Proposition \ref{lemma-compl-P1-bound-s}
the log divisor $K+B$ is 
$n$-complemented for some 
$n \in \NNN_{d}(\R)$.
For this complement $B^+$, we have $B^+\ge B$.
Since $K_{X}+B\equiv 0$, $B^+=B$. In particular,
$nB$ is integral.
Thus we can put $\lambda:=1/\lcm (\NNN_{d-1}(\R))$.
This proves the statement in case (ii).
\end{proof}

\begin{corollary}
\label{cor-Adj-hor-mult-1}
Notation and assumptions as in Lemma \xref{lemma-Adj-hor-mult-1}.
Let $E$ be any divisor \textup(exceptional or not\textup) 
over $X$. Then either $a(E,X,B)=1$ or 
$a(E,X,B)\le 1-\lambda$. 
\end{corollary}
\begin{proof}
We can take $\lambda<\ep'$.
If $1-\lambda<a(E,X,B)<1$, replace 
$(X,B)$ with a crepant blowup of $E$, see \cite[21.6.1]{Utah} and apply 
Lemma \xref{lemma-Adj-hor-mult-1}.
\end{proof}

\subsection{Proof of Theorem \ref{main-result0}
in the case when $(X,B)$ is not klt}
\label{subsect-1}
Let $(X,B)$ be a $0$-pair such that $X$ is FT and
$B\in \Phi(\R)$.
We assume that $(X,B)$ is not klt.
Replace $(X,B)$ with its $\QQ$-factorial dlt modification.
Then in particular, $X$ is klt.
Moreover, $\down B\neq 0$ (and $B\in \Phi(\R)$).
Let $\lambda(d,\R)$ be as in Lemma \ref{lemma-Adj-hor-mult-1}
and let $0<\lambda<\lambda(d,\R)$.
If $X$ is not $\lambda$-lt, then for each exceptional divisor $E$ 
of discrepancy $a(E,X,0)< -1+\lambda$ 
by Corollary \ref{cor-Adj-hor-mult-1}
we have $a(E,X,B)=1$. Hence as in \ref{sub-reduction-1}
replacing $(X,B)$ with 
blowup of all such divisors $E$, see \cite[21.6.1]{Utah},
we get that $X$ is $\lambda$-lt and $(X,B)$
is a $0$-pair with $B\in \Phi(\R)$ and $\down{B}\neq 0$. 
Run $K$-MMP:\quad $X\dashrightarrow X'$ and let $B'$ 
be the birational transform of $B$.
Since $B\neq 0$, $X'$ admits a $K$-negative Fano fibration $X'\to Z'$
over a lower-dimensional variety $Z'$.
By our construction $X'$ is $\lambda$-lt and 
$(X',B')$
is a non-klt $0$-pair with $B'\in \Phi(\R)$. By Proposition 
\ref{cor-pull-back_compl-I} we can pull-back $n$-complements 
from $X'$ to $X$ if $I(\R) \mid n$. 
If $Z'$ is a point, then 
$\rho(X')=1$ and $X'$ is a klt Fano variety. In this case,
arguing as in 
\ref{para-rem-Noetherian}
we get that $(X',\Supp B')$ belong to an algebraic family.
By Proposition \ref{prop-1} and Noetherian induction 
the log divisor $K+B$ has a bounded $nI$-complement
for some $n\le C(d,\R)$. 
Finally, if $\dim Z'>0$, then we apply 
Proposition \ref{prop-D-not-big} below.
This will be explained in \ref{subsetion-0-pairs-proof}.
Theorem \ref{main-result0}
in the non-klt case is proved.

\section{Approximation and complements}
\label{approximation}
The following Lemma \xref{open} shows that the
existence of $n$-complements is an open condition in the space of
all boundaries $B$ with fixed $\Supp{B}$.
\subsection{Notation}
Let $\BBB$ be a finite set of prime divisors $B_i$.
Recall that $\mathfrak{D}_{\BBB}$ denotes
the $\RR$-vector space all $\RR$-Weil divisors $B$
with $\Supp B=\sum_{B_i\in \BBB} B_i$, where the $B_i$ are prime divisors. 
Let
\[
\mathfrak{I}_{\BBB}:=\left\{\sum
\beta_iB_i\in\mathfrak{D}_{\BBB}\mid 0\le\beta_i\le 1,\ \forall
i\right\}
\]
be the unit cube in $\mathfrak{D}_{\BBB}$.

\begin{lemma}
\label{open}
Let $(X,B)$ be a log pair where $B$ is an
$\RR$-boundary. Assume that $K+B$ is $n$-complemented. Then there
is a constant $\ep=\ep(X,B,n)>0$ such that $K+B'$ is also
$n$-complemented for any $\RR$-boundary $B'\in \mathfrak{D}_B$
with $\|B-B'\|<\ep$.
\end{lemma}
\begin{proof}
Let $B^+=B^{\sharp}+\Lambda$ be an $n$-complement, where $\Lambda$
and $B$ have no common components and
$B^{\sharp}\in\mathfrak{D}_B$. Write $B=\sum b_iB_i$, $B'=\sum
b_i'B_i$, $B^{\sharp}=\sum b_i^+B_i$. Take $\ep$ so that
\[
0<(n+1)\ep<\min(1-\fr{(n+1)b_i} \mid 1\le i\le r,\quad b_i<1).
\]
We claim that $B^+$ is also an $n$-complement of $B'$ whenever
$\|B-B'\|<\ep$. If $b_i=1$, then obviously $b_i^+=1$. So, it is
sufficient to verify the inequalities $nb_i^+\ge \down{(n+1)b_i'}$
whenever $b_i^+<1$ and $b_i<1$. Indeed, in this case,
\[
\down{(n+1)b_i'}\le
\down{(n+1)b_i+(n+1)(b_i'-b_i)}=\down{(n+1)b_i}\le nb_i^+.
\]
(because $(n+1)(b_i'-b_i)<(n+1)\ep<\min(1-\fr{(n+1)b_i})$). This
proves the assertion.
\end{proof}

\begin{corollary}
\label{corollary-open}
For any $D\in Z_{d-1}(X)$, the subset
\[
\mathfrak{U}^n_D:=\{B\in \mathfrak{I}_D\mid \text{$K+B$ is
$n$-complemented}\}
\]
is open in $\mathfrak{I}_D$.
\end{corollary}

\begin{proposition}
\label{prop-1}
Fix a positive integer $I$.
Let $X$ be an FT variety such that $K_X$ is $\QQ$-Cartier 
and let 
$B_1,\dots,B_r$ are $\QQ$-Cartier divisors on $X$.
Let $\BBB:=\sum_{i=1}^r B_i$.
Then for any boundary $B\in\mathfrak{I}_{\BBB}$ such that $K+B$ is
lc and $-(K+B)$ is nef, there is an $n$-complement of $K+B$ for
some $n\le\Const\left(X,\BBB\right)$ and $I\mid n$.
\end{proposition}

\begin{lemma}
\label{rem-prop-1}
In notation of Proposition \xref{prop-1} the following holds.
\para
\label{condition} {\rm (Effective base point freeness)}
There is a positive integer $N$ such that for any integral nef Weil
$\QQ$-Cartier divisor of the form $mK+\sum m_iB_i$ the linear
system $\left|N(mK+\sum m_iB_i)\right|$ is base point free.
\end{lemma}

\begin{proof}
Indeed, we have $\Pic(X)\simeq \ZZ^{\rho}$ (see, e.g.,
\cite[Prop. 2.1.2]{Iskovskikh-Prokhorov-1999}).
In the space $\Pic(X)\otimes \RR\simeq \RR^{\rho}$ we
have a closed convex cone $\operatorname{NEF}(X)$, the cone of nef divisors.
This cone is dual to the Mori cone $\NE(X)$, so it is rational polyhedral
and generated by a finite number of semiample Cartier divisors
$M_1,\dots, M_s$. Take a positive integer $N'$ so that all the linear systems
$|N'M_i|$ are base point free, and $N'K$, $N'B_1,\dots, N'B_r$ are
Cartier. Write
\[
N'K\q{\QQ} \sum_{i=1}^s \alpha_{i,0}M_i,\quad N'B_j\q{\QQ} \sum_{i=1}^s
\alpha_{i,j}M_i,\qquad \alpha_{i,j}\in\QQ,\qquad \alpha_{i,j}\ge
0.
\]
Let $N''$ be the common multiple of denominators of the
$\alpha_{i,j}$. Then
\begin{multline*}
{N'}^2N''\left(mK+\sum_{j=1}^r m_jB_j\right)\sim
N'N''m\left(\sum_{i=1}^s \alpha_{i,0}M_i\right)+\\ \sum_{j=1}^r
N'N''m_j \left(\sum_{i=1}^s \alpha_{i,j}M_i\right)= \sum_{i=1}^s
\left(mN''\alpha_{i,0}+\sum_{j=1}^r N''m_j
\alpha_{i,j}\right)N'M_i
\end{multline*}
The last (integral) divisor generates a base point free linear
system, so we can take $N={N'}^2N''$.
\end{proof}

In the proof of Proposition \xref{prop-1} we follow arguments of
\cite[Example 1.11]{Shokurov-2000}, see also \cite[5.2]{Shokurov-1992-e}.
\begin{proof}[Proof of Proposition {\xref{prop-1}}]
Define the set
\begin{equation}
\label{def-MMM}
\MMM=\MMM_{\BBB}:=\left\{B\in \mathfrak{I}_{\bar
B}\quad\left|\quad \text{$K+B$ is lc and $-\left(K+B\right)$ is
nef}\right.\right\}
\end{equation}
Then $\MMM$ is a closed compact convex polyhedron in
$\mathfrak{I}_{\BBB}$. It is sufficient to show the existence of
some $n$-complement for any $B\in \MMM$. Indeed, then $\MMM\subset
\bigcup\limits_{n\in\ZZ_{>0}}\mathfrak{U}^n_{\BBB}$. By taking a finite
subcovering $\MMM\subset
\bigcup\limits_{n\in\mathcal{S}}\mathfrak{U}^n_{\BBB}$ we get a
finite number of such $n$.

Assume that there is a boundary $B^o=\sum_{i=1}^r b_i^oB_i\in
\MMM$ which has no any complements. By \cite[Ch. 1, Th. VII]{Cassels-1957}
there is infinite many rational points $(m_1/q,\dots,m_r/q)$ such
that
\[
\max\left(\left|\frac{m_1}{q}-b_1^o\right|, \dots,
\left|\frac{m_r}{q}-b_r^o\right|\right)
<\frac{r}{(r+1)q^{1+1/r}}<\frac1{q^{1+1/r}}.
\]
Denote $b_i:=m_i/q$ and $B:=\sum b_iB_i$. Thus,
$\|B-B^o\|<1/q^{1+1/r}$. Then our proposition is an easy
consequence of the following

\begin{claim}
For $q\gg 0$ one has
\para
\label{claim-11}
$\down{(qN+1)b_i^o}\le qNb_i$ whenever $b_i<1$\textup;
\para\label{claim-2}
$B\equiv B^o$ and $-(K+B)$ is nef\textup; and
\para
\label{claim-3}
$K+B$ is lc.
\end{claim}

Indeed, by \eqref{condition} the linear system $|-qN(K+B)|$ is
base point free. Let $F\in |-qN(K+B)|$ be a general member. Then
$K+B+\frac1{qN}F$ is an $qN$-complement of $K+B^o$, a
contradiction.
\end{proof}
\begin{proof}[Proof of Claim]
By the construction
\[
\down{(qN+1)b_i^o}=m_iN+\down{b_i^o+qN(b_i^o-b_i)}
\]
Put $c:=\max\limits_{b_i^o<1}\{ b_1^o,\dots b_r^o\}$. Then for
$b_i<1$ we have $b_i^o<c<1$ and for $q\gg 0$,
\[
b_i^o+qN(b_i^o-b_i)<c+\frac{qN}{q^{1+1/r}}<1.
\]
This proves \eqref{claim-11}.

Further, let $L_1,\dots, L_r$ be a finite set of curves generating
$N_1(X)$. We have
\begin{multline*}
\left|L_j\cdot(B-B^o)\right|= \left|\sum_i\frac{m_i}q(L_j\cdot
B_i)- \sum_i b_i^o(L_j\cdot B_i)\right| \\ < \frac1{q^{1+1/r}}
\sum_i(L_j\cdot B_i).
\end{multline*}
If $q\gg 0$, then the right hand side is $\ll 1/q$ while the left
hand side is from the discrete set $\pm L_j\cdot
(K+B^o)+\frac{1}{qN}\ZZ$ (because $qB$ is an integral divisor and by
our assumption \eqref{condition}). Hence the left hand side is zero
and $B\equiv B^o$. This proves \eqref{claim-2}.

Finally, we have to show that $K+B$ is lc. Assume the converse. By
\eqref{condition} the divisor $qN\left(K+B\right)$ is Cartier. So
there is a divisor $E$ of the field $\Bbbk(X)$ such that $a(E,X,B)\le
-1-1/qN$ and $a(E,X,B^o)\ge -1$. On the other hand, $a(E,X,\sum
\beta_iB_i)$ is an affine linear function in $\beta_i$:
\[
\frac{1}{qN}\le a(E,X,B^o)-a(E,X,B)=\sum
c_i(b_i^o-b_i)<\frac{\Const}{q^{1+1/r}}
\]
which is a contradiction.
\end{proof}

The following is the first induction step to prove 
Theorem \ref {main-result}.

\begin{corollary}[One-dimensional case]
\label{lemma-compl-P1-bound}
Fix a finite set $\R\subset [0,\, 1]\cap \QQ$ and a positive integer 
$I$.
Then the set $\NNN_1(\R)$ is finite.
\end{corollary}

\begin{proof}
Let $(X,B)$ be a one-dimensional log pair satisfying conditions 
of \eqref{para-def-N}. Since $X$ is FT, $X\simeq\PP^1$. 
Since $B\in \Phi(\R)$ and $\R$ is finite, 
we can write $B=\sum_{i=1}^r b_iB_i$, where $b_i\ge \delta$
for some fixed $\delta>0$.
Thus we may assume that $r$
is fixed and $B_1,\dots,B_r$ are fixed distinct points.
Then by Proposition \ref{prop-1} we have a desired complements.
\end{proof}

\begin{example}
Let $X\simeq \PP^1$.
If $\R=\{0,\, 1\}$, then $I(\R)=1$ and $\Phi(\R)$
is the set of standard multiplicities. In this case,
it is easy to compute that
$\NNN_1(\R)=\{1,\,2,\,3,\,4,\,6\}$
\cite[5.2]{Shokurov-1992-e}.
Consider more complicated case when
$\R=\{0,\, \frac12,\, \frac 23,\, \frac 34,\, \frac 56,\,
1\}$. Then $I=12$ and one can compute that
\[
\NNN_1(\R)=12\cdot \{1,\, 2,\, 3,\, 4,\, 5,\, 7,\, 8,\, 9,\, 11\}.
\]
Indeed, assume that $(X,D)$ has no any $12n$-complements for 
$n\in $ $\{1,  $ $2,$  $3, $ $4, $  $ 5, $ $7,$ $ 8,$ $ 9,$ $ 11\}$.
Write $D=\sum_{i=1}^r d_iD_i$, where $D_i\neq D_j$ for $i\neq j$.
It is clear that the statement about the existence of
an $n$-complement $D^+$ such that $D^+\ge D$ is equivalent to
the following inequality
\begin{equation}
\label{eq-compl-P1}
\sum_i \up{nd_i}\le 2n.
\end{equation}
Since $d_i=1-r_i/m_i$, where
$r_i\in \R$, $m_i\in \ZZ_{>0}$, we have $d_i\ge 1/6$ for all $i$.
We claim that at least one denominator of
$d_i$ does not divide $24$. Indeed,
otherwise $24D$ is an integral divisor and $D^+:=D+\frac 1{24}\sum_{j=1}^k D_j$
is an $24$-complement, where $D_j\in X$ are general points and
$k=24(2-\deg D)$.
Thus we may assume that the denominator of $d_1$ does
not divide $24$.
Since $d_1=1-r_1/m_1$, where $r_1\in \R$, we have
$m_1\ge 3$ and the equality holds only if $r_1=2/3$ or $5/6$.
In either case, $d_1\ge 13/18$.

Recall that a log pair $(X,D)$ of global type is said to be \textit{exceptional}
if at has at least one $\QQ$-complement and
any $\QQ$-complement is klt.
If $(X,D)$ is not exceptional, we can increase $d_1$ by putting $d_1=1$.
Then as above $13/18\le d_2\le 5/6$, so $r=3$ and $d_3\le 5/18$.
Now there are only a few possibilities for $d_2$ and $d_3$:
\begin{center}
\begin{tabular}{c|ccccccc}
$d_2$&$\frac{13}{18}$ &$\frac34$ &$\frac79$
&$\frac{19}{24}$ &$\frac45$ &$\frac{13}{16}$&$\frac56$
\\[8pt]
\hline
\\
$d_3$ &$\le\frac5{18}$ &$\le\frac14$ &$\le\frac29$ &$\le\frac5{24}$
&$\le\frac15$ &$\le\frac3{16}$ &$\le\frac16$
\\
\end{tabular}
\end{center}
\par\medskip\noindent
In all cases $K+D$ has a $12n$ complement for some $n\in \{1,\, 2,\, 3,\, 4,\, 5\}$.
In the exceptional case, there is a finite number of
possibilities for $(d_1,\dots,d_r)$.
However the computations are much longer.
We omit them.
\end{example}

\section{The main theorem: Case $\rho=1$}
\label{sec-rho=1}
\subsection{}
\label{rho=1-begin}
Now we begin to consider case
\eqref{case-rho=1}. Thus we assume that $(X,B)$ 
is klt but not
$\ep'$-lt, where we take $\ep'$ so that 
$\ep'N=1$ for some integer $N\ge N_{d-1}(\ov \R)+2$.
Then obviously, $\ep_{d-1}(\ov \R)\ge \ep'>0$. Further, assume that 
applying the general reduction from Section \ref{sect-reduction}
we get a pair $(Y,B_Y)$, where $Y$ is FT and $\QQ$-factorial,
$\rho(Y)=1$, the $\QQ$-divisor $-(K_Y+B_Y)$ is nef, and 
\[
\discr(Y,B_Y)>-1+\ep'.
\]
Moreover,
$K_Y+D_Y$ is ample and $\down{D_Y}\neq 0$, where $D_Y$ is a boundary 
satisfying \eqref{eq-D}.
In particular, $B_Y\neq 0$.

\subsection{}
Assume that the statement of Theorems \ref{main-result}
and \ref{main-result0}
is false in this case. Then there is a sequence of klt log pairs
$(X^{(m)},B^{(m)})$ as in \ref{rho=1-begin}
and such that complements of
$K_{X^{(m)}}+B^{(m)}$ are unbounded.
More precisely, for each $K_{X^{(m)}}+B^{(m)}$, let 
$n_m$ be the minimal positive integer such that $I(\R)\mid n_m$ 
and $K_{X^{(m)}}+B^{(m)}$ is 
$n_m$-complemented. We assume that the sequence $n_m$
is unbounded.
We will derive a contradiction.

\subsection{}
By Corollary \ref{lemma-compl-P1-bound} we have $\dim X^{(m)}\ge 2$.
By our hypothesis we have 
a sequence of birational maps $X^{(m)}\dashrightarrow Y^{(m)}$,
where $Y^{(m)}$, $X^{(m)}$, $B^{(m)}$ and $D^{(m)}$ are as above 
for all $m$.
Recall that by Corollary \ref{cor-reduct-discrs}
\begin{equation}
\label{eq-rho=1-discr}
\discr(Y^{(m)})\ge \discr(Y^{(m)},B^{(m)}_Y)
\ge\discr (X^{(m)},B^{(m)})> -1+\ep',
\end{equation}
where $0<\ep' \le \ep_{d-1}(\ov \R)$.
Note that $-(K_{Y^{(m)}}+B^{(m)})$ is nef, so
by Conjecture \ref{BAB} the sequence of varieties $Y^{(m)}$ is bounded.
By Noetherian induction (cf. \eqref{para-rem-Noetherian}),
we may assume that $Y^{(m)}$ is fixed, that is, $Y^{(m)}=Y$.

Let $Y\hookrightarrow\PP^N$ be an embedding and let
$H$ be a hyperplane section of $Y$. Note that
the multiplicities of $B_Y^{(m)}=\sum b_i^{(m)}B_i^{(m)}$ are bounded from below:
$b_i^{(m)}\ge \delta_0>0$, where $\delta_0:=\min \Phi(\R)\setminus \{0\}$.
Then, for each $B_i^{(m)}$,
\[
\delta_0 H^{d-1}\cdot B_i^{(m)}\le H^{d-1}\cdot B^{(m)}_Y\le -H^{d-1}\cdot K_Y.
\]
This shows that the degree of $B_i^{(m)}$ is bounded
and $B_i^{(m)}$ belongs to an algebraic family. Therefore we may assume
that $\Supp{B_Y^{(m)}}$ is also fixed: $B_i^{(m)}=B_i$.

\subsection{}
\label{subsect-rho=1-case1}
Assume that the multiplicities of $B^{(m)}_Y$ are
bounded from $1$, i.e., $b_i^{(m)}\le 1-c$, where
$c>0$. Then we
argue as in \xref{sufficient-condition}.
By Claim
\xref{red-not-contr} divisorial contractions in \eqref{MMP} do not
contract components of $B^{(m)}$ with multiplicities $b_i^{(m)}\ge
1-\ep'$. Therefore the multiplicities of $B^{(m)}$ on $X^{(m)}$
are also bounded
from the $1$. Combining this with $\discr (X,B)>-1+\ep'$
and Conjecture \xref{BAB} we get that $(X,\Supp B)$ belong to an
algebraic family. By Noetherian induction (cf. \eqref{para-rem-Noetherian}) 
we may assume that
$(X,\Supp B)$ is fixed. Finally, by Proposition \xref{prop-1} 
we have that $(X,B)$ has a bounded
$n$-complement such that $I(\R)\mid n$.

Thus by our construction the only possibility is the case below.

\subsection{}
From now on we consider the remaining case when
some multiplicity of $B^{(m)}_Y$ is accumulated to $1$
and we will derive a contradiction.
Since $\Supp B^{(m)}_Y$ does not depend on $m$,
by passing to a subsequence we may assume that
the limit $B^{\infty}_Y:=\lim\limits_{m\to\infty} B^{(m)}_Y$ exists and
$\down{B^{\infty}_Y}\neq 0$.
As above, write $B^{\infty}_Y=\sum b^{\infty}_i B_i$.
Up to permutations of components we may assume that 
$b^{\infty}_1=1$.
It is clear that
$-(K_Y+B^{\infty}_Y)$ is nef and $\discr(Y,B^{\infty}_Y)\ge -1+\ep'$.
In particular, $(Y,B^{\infty}_Y)$ is plt.

\begin{claim}
\label{claim-bi}
Under the above hypothesis, we have $b_j^{\infty}\le 1-\ep'$ for
all $1< j\le r$. Moreover,
by passing to a subsequence we may assume the following\textup:
\begin{enumerate}
\item
If $b_j^{\infty}=1$, then $j=1$ and $b_1^{(m)}$ is strictly increasing.
\item
If $b_j^{\infty}<1-\ep'$, then $b_j^{(m)}=b_j^{\infty}$ is a constant.
\item
If $b_j^{\infty}=1-\ep'$, then $b_j^{(m)}$ is either a constant
or strictly decreasing.
\end{enumerate}
In particular, $B^{\infty}_Y\in\Phi(\R)$, 
$B^{\infty}_Y$ is a $\QQ$-boundary, and $D_Y^{(m)}\ge B^{\infty}_Y$ for $m\gg 0$.
\end{claim}
\begin{proof}
Since $\rho(Y) =1$ and $Y$ is $\QQ$-factorial, the
intersection $B_1\cap B_j$ on $Y$ is of codimension two and non-empty. For a general
hyperplane section $Y\cap H$, by \eqref{eq-rho=1-discr} we have
the inequality
\[
\discr(Y\cap H,B^{\infty}\cap H)\ge -1+\ep'.
\]
Thus by Lemma
\xref{pair-discr} below, we have $b_1^{\infty}+b_j^{\infty}\le 2-\ep'$,
i.e., $b_j^{\infty}\le 1-\ep'$ for
all $1< j\le r$. The rest follows from the fact that the set $\Phi(\R)\cap [0,1-\ep']$ is finite.
\end{proof}

\begin{lemma}[cf. {\cite[Prop. 5.2]{Shokurov-2000}, \cite[\S 9]{Prokhorov-2001}}]
\label{pair-discr}
Let $(S\ni
o,\Lambda=\sum \lambda_i\Lambda_i)$ be a log surface
germ. Assume that
$\discr(S,\Lambda)\ge -1+\ep$ at $o$
for some positive $\ep$. Then $\sum \lambda_i\le 2-\ep$.
\end{lemma}
\begin{proof}
Locally near $o$ there is an \'etale outside of
$o$ Galois cover $\pi\colon S'\to S$ such that $S'$ is smooth.
Let $\Lambda':=\pi^*\Lambda$ and $o':=\pi^{-1}(o)$. Then 
$\discr(S',\Lambda')\ge\discr(S,\Lambda)\ge -1+\ep$ at $o$
(see, e.g., \cite[Proposition 20.3]{Utah}).
Consider the blow up of $o'\in S'$. 
We get an exceptional divisor $E$ of
discrepancy
\[
-1+\ep\le
a(E,S',\Lambda')=1-\sum \lambda_i.
\]
This gives us the desired inequality.
\end{proof}

\begin{corollary}
$b_j^{\infty}=1-\ep'$ for some $j$.
\end{corollary}
\begin{proof}
Indeed, otherwise $D^{(m)}_Y=B^{\infty}_Y$ for $m\gg 0$ and
$-(K_{Y}+D^{(m)}_Y)$ is nef, a contradiction.
\end{proof}

We claim that the log divisor 
$K_{Y}+B^{\infty}_Y$ is $n$-complemented, where
$n\in \NNN_{d-1}(\ov \R)$.
Recall that $b_1^{\infty}=1$.
Put $B':=B^{\infty}-B_1$. By the last corollary
$B'\neq 0$. By Proposition \ref{prop-prop}
we have $\Diff_{B_1}(B')\in \Phi(\ov \R)$.
Recall that $-(K_Y+B^\infty)$ is nef.
Since 
$(Y,B^{\infty}_Y)$ is plt, the pair $(B_1,\Diff_{B_1}(B'))$ 
is klt. 
Further, $-(K_{B_1}+\Diff_{B_1}(0))$ is ample (because $\rho(Y)=1$), so
$B_1$ is FT.
Thus by the inductive hypothesis
there is an
$n$-complement $K_{B_1}+\Diff_{B_1}(B')^+$ of $K_{B_1}+\Diff_{B_1}(B')$
for some $n\in \NNN_{d-1}(\ov \R)$.
By Lemma \ref{lemma-PPP-n} $B'\in \PPP_n$.
Take a sufficiently small positive $\delta$ and let $j$ be such that
$b_j^{\infty}=1-\ep'$. We claim that
$B'-\delta B_j\in \PPP_n$. Indeed, otherwise $(n+1)b_j^\infty$
is an integer. On the other hand,
\[
n+1> (n+1)b_j^\infty=(n+1)(1-\ep')\ge n+1-(n+1)/N.
\]
where $N\ge N_{d-1}(\ov \R)+2\ge n+2$. This is impossible.
Thus, $B'-\delta B_j\in \PPP_n$. Since $-(K_Y+B^\infty-\delta B_j)$ is
ample, by Proposition \ref{prodolj} 
the $n$-complement $K_{B_1}+\Diff_{B_1}(B')^+$ of
$K_{B_1}+\Diff_{B_1}(B'-\delta B_j)$ can be extended to $Y$.
So there is an $n$-complement $K_Y+B^{+}$ of 
$K_Y+B^{\infty}-\delta B_j$.
Write $B^+=\sum b^+_i B_i$.
Since $B-\delta B_j\in \PPP_n$, we have $B^{+}\ge B^{\infty}-\delta B_j$.
Moreover, since $nb_j^+$
is an integer and $1\gg \delta>0$, we have $b_j^+\ge b_j^\infty$.
Hence $B^+\ge B^{\infty}_Y$ and $B^+$ is also an
$n$-complement of $K_Y+B^{\infty}$ (see
Remark
\ref{rem-complne}).
By Lemma \ref{open} \ $K_Y+B^{+}$ is also an
$n$-complement of $K_Y+B^{(m)}$ for $m\gg 0$.

By Lemmas \ref{lemma-PPP-n-1}
and \ref{lemma-PPP-n} we have $B^+\ge B^{(m)}$
for $m\gg 0$. More precisely,
\[
b_i^+\quad
\begin{cases}
\quad=1& \text{if $b_i^{\infty}\ge 1-\ep'$,}
\\
\quad\ge b_i^{(m)}=b_i^{\infty}& \text{if $b_i^{\infty}< 1-\ep'$.}
\end{cases}
\]
By the construction of $D$ we have $D^{(m)}\le B^+$. 
Hence $-(K_Y+D^{(m)})$ is nef, a contradiction.
This completes the proof of Theorems \ref{main-result}
and \ref{main-result0} in 
case \eqref{case-rho=1}.

\section{Effective adjunction}
\label{sect-Adj}
\renewcommand{\thefootnote}{\Roman{footnote}}
In this section we discuss 
the adjunction conjecture for fibre spaces.
This conjecture can be considered as a generalization 
of the classical Kodaira canonical bundle formula 
for canonical bundle, see 
\cite{Kodaira-1963}, 
\cite{Fujita-1986ZE},
\cite{Kawamata-1997-Adj},
\cite{Kawamata-1998},
\cite{Ambro-PhD},
\cite{Fujino-1999app},
\cite{Fujino-Mori-2000},
\cite{Fujino-2003-AG},
\cite{Ambro-2004S}, 
\cite{Ambro-2005}.

\subsection{The set-up}
\label{not-adj}
Let $f\colon X\to Z$ be a surjective 
morphism of normal varieties and let
$D=\sum d_iD_i$ be an $\RR$-divisor on $X$ such that $(X,D)$ is lc
near the generic fibre of $f$
and $K+D$ is $\RR$-Cartier over 
the generic point of any prime divisor $W\subset Z$. 
In particular, $d_i\le 1$ whenever
$f(D_i)=Z$. Let $d:=\dim X$ and $d':=\dim Z$.

For any divisor $F=\sum \alpha_i F_i$ on $X$,
we decompose $F$ as $F=F^\h+F^\ver$, where 
\[
F^\h:=\sum\limits_{f(F_i)=Z} \alpha_i F_i, \qquad
F^\ver:=\sum\limits_{f(F_i)\neq Z} \alpha_i F_i. 
\]
These divisors $F^\h$ and $F^\ver$ are called
the \textit{horizontal} and \textit{vertical} parts of $F$, respectively.

\subsection{Construction}
\label{def-cW}
For a prime divisor $W\subset Z$, define a real number $c_W$ as the log
canonical threshold over the generic point of $W$:
\begin{equation}
\label{cW}
c_W:=\sup\left\{c \mid (X,D+cf^\bullet W)\quad \text{is lc over the
generic point of $W$}\right\}.
\end{equation}
It is clear that $c_W\in\QQ$ whenever $D$ is a $\QQ$-divisor. Put
$d_W:=1-c_W$. Then the $\RR$-divisor
\[
D\sdiv:=\sum_W d_WW
\]
is called the \emph{divisorial part}
of adjunction (or \emph{discriminant} of $f$) for
$K_X+D$. 
It is easy to see that $D\sdiv$ is a divisor, i.e., $d_W$ is zero 
except for a finite number of prime divisors.

\begin{remark}
\label{rem-adj-2}
\begin{enumerate}
\item 
Note that the definition of the discriminant $D\sdiv$ is a codimension
one construction, so computing $D\sdiv$ we can systematically remove
codimension two subvarieties in $Z$ and pass to general hyperplane
sections $f_H\colon X\cap f^{-1}(H)\to Z\cap H$.
\item
Let $h\colon X'\to X$ be a birational contraction and let $D'$ be
the crepant pull-back of $D$:
\[
K_{X'}+D'=h^*(K_X+D),\quad h_*D'=D.
\]
Then ${D'}\sdiv=D\sdiv$, i.e., the discriminant $D\sdiv$ does not depend on the
choice of crepant birational model of $(X,D)$ over $Z$.
\end{enumerate}
\end{remark}

The following lemma is an immediate consequence of the definition.

\begin{lemma}
\label{lemma-adj-prelim}
Notation as in \xref{not-adj}.
\begin{enumerate}
\item {\rm (effectivity, cf. \cite[3.2]{Shokurov-1992-e})}
If $D$ is boundary over the generic point 
of any prime divisor $W\subset Z$, then $D\sdiv$ effective.
\item {\rm (semiadditivity, cf. \cite[3.2]{Shokurov-1992-e})}
Let $\Delta$ be an $\RR$-divisor on $Z$ and let $D':=D+f^\bullet\Delta$.
Then ${D'}\sdiv=D\sdiv+\Delta$.
\item
$(X,D)$ is klt \textup(resp., lc\textup) over the generic point of $W$ if and only if
$d_W< 1$ \textup(resp., $d_W\le 1$\textup).
\item
If $(X,D)$ is lc and $D$ is an $\RR$- \textup(resp., $\QQ$-\textup)boundary, then
$D\sdiv$ is an $\RR$- \textup(resp., $\QQ$-\textup)boundary.
\end{enumerate}
\end{lemma}

\subsection{Construction}
\label{construction-adj-def-mod}
From now on assume that $f$ is a contraction,
$K_X+D$ is $\RR$-Cartier, 
and $K+D\q{\RR} f^* L$
for some $\RR$-Cartier divisor $L$ on $Z$.
Recall that the latter means that there are real numbers $\alpha_j$ and
rational functions $\var_j\in \Bbbk(X)$ such that
\begin{equation}
\label{eq-def-dm}
K+D-f^*L = \sum \alpha_j\, (\var_j).
\end{equation}
Define the \textit{moduli part} $D\smod$ of
$K_X+D$ by 
\begin{equation}
\label{eq-def-dm-3}
D\smod:=L-K_Z-D\sdiv.
\end{equation}
Then we have 
\begin{equation}
\label{eq-adj-first-p-R}
K_X+D=f^*(K_Z+D\sdiv+D\smod)+\sum \alpha_j\, (\var_j).
\end{equation}
In particular,
\[
K_X+D\q{\RR} f^*(K_Z+D\sdiv+D\smod).
\]
Clearly,
$D\smod$ depends on the choice of representatives of 
$K_X$ and $K_Z$, and also on the choice of 
$\alpha_j$ and $\var_j$ in \eqref{eq-def-dm}.
Any change of $K_X$ and $K_Z$
and change of $\alpha_j$ and $\var_j$
gives a new $D\smod$ which differs from the original one
modulo $\RR$-linear equivalence.

If $K+D$ is $\QQ$-Cartier, the definition of 
the moduli part is more explicit. 
By our assumption \eqref{eq-def-dm}
there is a positive integer $I_0$ such that
$I_0(K+D)$ is linearly trivial
on the generic fibre. 
Then for some rational function $\psi\in \Bbbk(X)$, the divisor 
$M:=I_0(K+D)+(\psi)$ is vertical
(and $\QQ$-linearly trivial over $Z$). Thus,
\[
M-I_0f^*L=(\psi)+\sum I_0\alpha_j\, (\var_j),\qquad \alpha_j\in \QQ.
\]
Rewrite it in a more compact form:
$M-I_0f^*L=\alpha\, (\var)$,\ $\alpha\in \QQ$, $\var\in \Bbbk(X)$.
The function $\var$ vanishes on the generic fibre,
hence it is a pull-back of some function $\upsilon\in \Bbbk(Z)$.
Replacing $L$ with $L+\frac{\alpha}{I_0}\, (\upsilon)$ we get
$M=I_0f^*L$ and
\begin{equation}
\label{eq-def-dm-Q}
K_X+D-f^*L =\frac1{I_0}(\psi),
\qquad \psi\in \Bbbk(X).
\end{equation}
In other words, $K+D\q{I_0} f^*L$.
Here $L$ is $\QQ$-Cartier. 
Then again we define the moduli part $D\smod$ of
$K_X+D$ by \eqref{eq-def-dm-3}, where 
$L$ 
is taken to satisfy \eqref{eq-def-dm-Q}. In this case,
$D\smod$ is $\QQ$-Cartier and we have 
\begin{equation}
\label{eq-adj-first-p-Q}
K_X+D=f^*(K_Z+D\sdiv+D\smod)+\frac 1{I_0}\, (\psi).
\end{equation}
In particular, 
\[
K_X+D\q{I_0}f^*(K_Z+D\sdiv+D\smod).
\]
As above,
$D\smod$ depends on the choice of representatives of 
$K_X$ and $K_Z$, and also on the choice of 
$I_0$ and $\psi$ in \eqref{eq-def-dm-Q}.
Note that $I_0$ depends only on $f$ and the horizontal part 
of $D$. Once these are fixed, we usually will assume that
$I_0$ is a constant. Then any change of 
$K_X$, $K_Z$, and $\psi$
gives a new $D\smod$ which differs from the original one
modulo $I_0$-linear equivalence.

\begin{remark1}
By Lemma \ref{lemma-adj-prelim}
$(D+f^\bullet \Delta)\smod=D\smod$.
Roughly speaking this means that 
``the moduli part depends only on the horizontal 
part of $D$''.
\end{remark1}

For convenience of the
reader we recall definition of b-divisors and related notions, see
\cite{Iskovskikh-2003-b-div}
for details.
\begin{definition}
Let $X$ be a normal variety.
Consider an infinite linear combination
$\bD:=\sum_P d_P P$, where $d_P\in \RR$ and $P$ runs through
all discrete valuations $P$ of the function field. 
For any birational model
$Y$ of $X$ define the trace of $\bD$ on $Y$ as follows
$\bD_Y:=\sum\limits_{\codim_Y P=1} d_P P$.
A \textit{b-divisor} is a linear combination $\bD=\sum_P d_P P$
such that the trace $\bD_Y$ on each birational model
$Y$ of $X$ is an $\RR$-divisor, i.e., only a finite number of 
multiplicities of $\bD_Y$ are non-zero. In other words, 
a b-divisor is an element of
$\underleftarrow{\lim} \Div_{\RR} (Y)$, where 
$Y$ in the inverse limit runs through all normal
birational models $f\colon Y\to X$, $\Div_{\RR} (Y)$
is the group of $\RR$-divisors of $Y$, and 
the map
$\Div_{\RR} (Y)\to \Div_{\RR} (X)$ is the push-forward. 
Let $D$ be a $\RR$-Cartier divisor on $X$. 
The \textit{Cartier closure} of $D$ is a b-divisor $\ov D$ whose trace 
on every birational model $f\colon Y\to X$ is $f^*D$.
A b-divisor $\bD$ is said to be b-Cartier 
if there is a model $X'$ 
and a $\RR$-Cartier divisor $D'$ on $X'$ such that
$\bD=\ov {D'}$. 
A b-divisor $\bD$ is said to be b-nef 
(resp. b-semiample, b-free)
if it is b-Cartier and 
there is a model $X'$ 
and a $\RR$-Cartier divisor $D'$ on $X'$ such that
$\bD=\ov {D'}$ and $D'$ is nef (resp. semiample, integral and free).

$\QQ$- and $\ZZ$-versions of b-divisors are defined similarly.
\end{definition}

\begin{remark}
\label{rem-adj-3}
Let $g\colon Z'\to Z$ be a birational contraction. Consider the following diagram
\begin{equation}
\label{eq-Diff-sq-diag}
\begin{CD}
X'@>h>> X
\\
@V{f'}VV @V{f}VV
\\
Z'@>{g}>> Z
\end{CD}
\end{equation}
where $X'$ is a resolution of the dominant component of $X\times _{Z} Z'$.
Let $D'$ be the crepant pull-back of $D$ that is 
$K_{X'}+D'=h^*(K_X+D)$ and $h_*D'=D$.
By Remark \ref{rem-adj-2} we have $g_*D'\sdiv=D\sdiv$. Therefore,
the discriminant defines a b-divisor $\bD\sdiv$.

For a suitable choice of $K_X'$, we can write
\[
h^* (K_X+D)=K_{X'}+D',
\]

Now we fix the choice of $K$, $\alpha_j$ and $\var_j$ in
\eqref{eq-def-dm} (resp. $K$ and $\psi$ in \eqref{eq-def-dm-Q})
and induce them naturally to $X'$. 
Then $D\smod$ and $D'\smod$ are uniquely determined
and $g_*D'\smod=D\smod$. This defines a b-divisor $\bD\smod$.

We can write
\[
K_{Z'}+D'\sdiv+D'= g^*(K_{Z}+D\sdiv+D\smod)+E,
\]
where $E$ is $g$-exceptional.
Since
\[
h^* f^*(K_Z+D\sdiv+D\smod) \equiv K_{X'}+D' \equiv
f'^*(K_{Z'}+D'\sdiv+D'\smod),
\]
we have $E=0$ (see \cite[1.1]{Shokurov-1992-e}),
i.e., $g$ is $(K_Z+D\sdiv+D\smod)$-crepant:
\begin{equation}
\label{eq-Diff-sq}
K_{Z'}+D'\sdiv+D'\smod = g^* (K_Z+D\sdiv+D\smod).
\end{equation}
\end{remark}

Let us consider some examples.
\begin{example}
Assume that the contraction $f$ is birational.
Then by the ramification formula 
\cite[\S 2]{Shokurov-1992-e}, \cite[Prop. 20.3]{Utah}
and negativity lemma \cite[1.1]{Shokurov-1992-e} 
we have $D\sdiv=f_*D$, $K+D=f^*(K_Z+D\sdiv)$,
and $D\smod=0$.
\end{example}

\begin{example}
Let $X=Z\times \PP^1$ and let $f$ be the natural projection
to the first factor. Take
very ample divisors $H_1,\dots, H_4$ on $Z$. Let $C$ be a section
and let $D_i$ be a general member of the linear system
$|f^*H_i+C|$. Put $D:=\frac12\sum D_i$. Then $K_X+D$ is
$\QQ$-linearly trivial over $Z$. By Bertini's theorem $D+f^*P$ is
lc for any point $P\in Z$.
Hence $D\sdiv=0$. On the other hand,
\[
K_X+D=f^*K_Z-2C+\frac12f^*\sum H_i +2C= f^*\left(K_Z+\frac12\sum
H_i\right).
\]
This gives us that $D\smod\q{\QQ}\frac12\sum H_i$.
\end{example}

\begin{example}
Let $X$ be a hyperelliptic surface. Recall that it is constructed
as the quotient $X=(E\times C)/G$ of the product of two elliptic
curves by a finite group $G$ acting on $E$ and $C$ so that the
action of $G$ on $E$ is fixed point free and the action on $C$ has
fixed points. Let
\[
f\colon X=(E\times C)/G\to \PP^1=C/G
\]
be the projection. It is clear that degenerate fibres of $f$ can
be only of type $m\mathrm{I}_0$. Using the classification of such
possible actions (see, e.g., \cite[Ch. V, Sect. 5]{BPV-1984}) we obtain
the following cases:
\par\medskip\noindent
\begin{center}
\begin{tabular}{p{85pt}p{100pt}p{100pt}}
{\rm Type}&{\rm singular fibres}&
\multicolumn{1}{c}{$D\sdiv$}\vspace{5pt}\\ \hline
\\
$\mathrm{a)}$\quad ($2K_X\sim 0$)& $2\mathrm{I}_0$,
$2\mathrm{I}_0$, $2\mathrm{I}_0$, $2\mathrm{I}_0$&
$\frac12P_1+\frac12P_2+\frac12P_3+\frac12P_4$\\ $\mathrm{b)}$\quad
($3K_X\sim 0$)& $3\mathrm{I}_0$, $3\mathrm{I}_0$, $3\mathrm{I}_0$&
$\frac23P_1+\frac23P_2+\frac23P_3$\\ $\mathrm{c)}$\quad ($4K_X\sim
0$)& $2\mathrm{I}_0$, $4\mathrm{I}_0$, $4\mathrm{I}_0$&
$\frac12P_1+\frac34P_2+\frac34P_3$\\ $\mathrm{d)}$\quad ($6K_X\sim
0$)& $2\mathrm{I}_0$, $3\mathrm{I}_0$, $6\mathrm{I}_0$&
$\frac12P_1+\frac23P_2+\frac56P_3$\\
\end{tabular}
\end{center}
\par\medskip\noindent
In all cases the moduli part $D\smod$ is trivial.
\end{example}

\begin{assumption}
\label{assumpt-adj-*}
Under the notation of \ref{not-adj} and \ref{construction-adj-def-mod} 
assume additionally that
$D$ is a $\QQ$-divisor and
there is a $\QQ$-divisor $\Theta$ on $X$ such that
$K_X+\Theta$ is $\QQ$-linearly trivial over $Z$ and
$(F,(1-t)D|_F +t \Theta|_F)$ is a klt log pair for any $0<t\le 1$, where
$F$ is the generic fibre of $f$. In particular, $\Theta$ and $D$ are
$\QQ$-boundaries near the generic fibre.
In this case, both $D\sdiv$ and $D\smod$ are $\QQ$-divisors. 
\end{assumption}

The following result is very important.
\begin{theorem}[{\cite{Ambro-2004S}}]
\label{th-Ambro-m}
Notation and assumptions as in
\xref{not-adj} and \xref{construction-adj-def-mod}.
Assume additionally that $D$ is a $\QQ$-divisor, $D$ is 
effective near the generic fibre, and
$(X,D)$ is klt near the generic fibre.
Then we have.
\begin{itemize}
\item[(i)] 
The b-divisor $\bK+\bD\sdiv$ is b-Cartier. 
\item[(ii)] 
The b-divisor 
$\bD\smod$ is b-nef.
\end{itemize}
\end{theorem}

According (ii) of Theorem \ref{th-Ambro-m} 
the b-divisor $\bD\smod$ is b-nef for $D\ge 0$ and $(X,D)$ is klt near
the generic fibre 
(see also \cite[Th.~2]{Kawamata-1998}, \cite{Fujino-1999app}).
We expect more.

\begin{conjecture}
\label{conj-main-adj}
Let notation and assumptions be as in \xref{not-adj} and \xref{assumpt-adj-*}.
We have

\para{} \label{conj-main-adj-1}
\textup(Log Canonical Adjunction\textup)
$\bD\smod$ is b-semiample.

\para{} \label{conj-main-adj-21}
\textup(Particular Case of Effective Log Abundance Conjecture\textup)
Let $X_\eta$ be the generic fibre of $f$.
Then $I_0(K_{X_\eta}+ D_\eta)\sim 0$,
where $I_0$ depends only on $\dim X_\eta$ and the
multiplicities of $D^\h$.

\para {} \label{conj-main-adj-2}
\textup(Effective Adjunction\textup)
$\bD\smod$ is effectively b-semiample, that is,
there exists a positive integer $I$ depending
only on the dimension of $X$ and
the horizontal multiplicities of $D$
\textup(a finite set of rational numbers\textup)
such that $I\bD\smod$ is very b-semiample, that is,
$I\bD\smod=\ov M$, where $M$ is a base point free
divisor on some model $Z'/Z$.
\end{conjecture}
Note that by \eqref{eq-adj-first-p-Q} we may assume that
\begin{equation}
\label{eq-adj-(*)}
K+D\q{I} f^*(K_Z+D\sdiv+D\smod).
\end{equation}

\begin{remark}
\label{rem-Adj_conj-compl}
We expect that 
hypothesis 
in \xref{conj-main-adj} can be weakened as follows.
\begin{enumerate}
\item[{\eqref{conj-main-adj-1}}]
It is sufficient to assume that $K+D$ is lc near 
the generic fibre, 
the horizontal part $D^{\h}$ of $D$ 
is an $\RR$-boundary, 
$K+D$ is $\RR$-Cartier, and 
$K+D\equiv f^*L$. 
\item[{\eqref{conj-main-adj-21}}]
$D^{\h}$ is a $\QQ$-boundary and 
$K+D\equiv 0$ near the generic fibre.
\item[{\eqref{conj-main-adj-2}}]
$D^{\h}$ is a $\QQ$-boundary,
$K+D$ is $\RR$-Cartier, and $K+D\equiv f^*L$.
\end{enumerate}
This however is not needed for the proof of the main theorem.
\end{remark}

\begin{remark}
In the notation of \eqref{conj-main-adj-21} we have $K_{X_\eta}+D_\eta\q{\QQ} 0$,
where $X_\eta$ is the generic fibre of $f$.
Assume that 
\begin{enumerate}
\item 
$X$ is FT, and
\item 
LMMP and conjectures \xref{BAB} and \xref{conj-main-adj} hold in 
dimensions $\le \dim X-\dim Z$.
\end{enumerate}
Then the pair 
$(X_\eta, D_\eta)$ satisfies the assumptions 
of Theorem \xref{main-result0}
with $\R$ depending only on horizontal
multiplicities of $D$. Hence $I_0(K_{X_\eta}+ D_\eta)\sim 0$,
where $I_0$ depends only on $\dim X_\eta$ and 
horizontal multiplicities of $D$.
Thus \eqref{conj-main-adj-21} holds automatically
under additional assumptions (i)-(ii).
\end{remark}

\begin{example}[Kodaira formula \cite{Kodaira-1963}, 
\cite{Fujita-1986ZE}]
\label{ex-elliptic}
Let $f\colon X\to Z$ be a fibration
satisfying \xref{construction-adj-def-mod}
whose generic fibre is an elliptic curve.
Then $D^\h=0$, $D=D^\ver$, and $I_0=1$. Thus we can write $K_X+D=f^*L$.
The $j$-invariant defines a rational map 
$J\colon Z \dashrightarrow \CC$.
By blowing up $Z$ and $X$ 
we may assume that both $X$ and $Z$ are smooth and $J$ is a morphism:
$J\colon Z\to \PP^1$.
Let $P$ be a divisor of degree $1$ on $\PP^1$.
Take a positive integer $n$ such
that $12n$ is divisible by the multiplicities of all the
degenerate fibres of $f$. In this situation, there is 
a generalization of the classical
Kodaira formula \cite{Fujita-1986ZE}:
\begin{equation*}
12n(K_X+D)=f^*\left(12nK_Z+12nD\sdiv+nJ^*P\right),
\end{equation*}
We can rewrite it as follows 
\begin{equation}
\label{eq-2}
K_X+D= f^*\left(K_Z+D\sdiv+\frac1{12}J^*P\right).
\end{equation}
Here $D\smod=\frac1{12}J^*P$ is semiample and the multiplicities of
$D\sdiv$ are taken from the table in Example~\xref{ex-Kodaira}
if $D=0$ over such divisors in $Z$ or $D$
is minimal as in Lemma \ref{lemma-adj-minimal}.
(Otherwise to compute $D\sdiv$ we can use
semiadditivity Lemma \ref{lemma-adj-prelim}, (ii).)
\end{example}

\begin{example}
Fix a positive integer $m$. 
Let $(E,0)$ be an elliptic curve with fixed group
low and let $e_m\in E$ be an $m$-torsion. Define the action of
$\muu_m:=\left\{ \sqrt[m]{1}\right\}$ on $E\times \PP^1$ by
\[
\ep(e,z)=(e+e_m,\ep z),\qquad e\in E,\ z\in \PP^1,
\]
where $\ep\in \muu_m$ is a primitive $m$-root. The quotient map
\[
X:=(E\times \PP^1)/\muu_m\longrightarrow \PP^1/\muu_m\simeq\PP^1
\]
is an elliptic fibration having exactly two fibres of types $m
\mathrm{I}_0$ over points $0$ and $\infty\in \PP^1$. Using the
Kodaira formula 
one can show that
\[
K_X=f^*K_Z+(m-1)F_0+(m-1)F_{\infty},
\]
where $F_0:=f^{-1}(0)_{\red}$ and
$F_{\infty}:=f^{-1}(\infty)_{\red}$. 
Hence, in \eqref{eq-def-dm-Q} we have $I_0=1$ and 
\[
L=K_Z+\left(1-\frac1m\right)\cdot 0+ \left(1-\frac1m\right)\cdot
\infty.
\]
Clearly,
\[
D\sdiv=\left(1-\frac1m\right)\cdot 0+ \left(1-\frac1m\right)\cdot
\infty.
\]
Hence $D\smod=0$, $I=I_0=1$, and $K_X=f^*(K_Z+D\sdiv)$.
\end{example}

\begin{corollary}[cf. {\cite[Th. 3.1]{Ambro-2004S}}]
\label{cor_Inv_aDj}
Let notation and assumptions be as in \xref{not-adj},
\xref{construction-adj-def-mod}, and \xref{assumpt-adj-*} 
\textup(cf. \xref{rem-Adj_conj-compl}\textup).
\begin{enumerate}
\item
If $(Z,D\sdiv+D\smod)$ is lc
and $\mathbf D\smod$ is effective, then 
$(X,D)$ is lc.
\item
Assume that \eqref{conj-main-adj-1} holds.
If $(X,D)$ is lc, then so is
$(Z,D\sdiv+D\smod)$ for a suitable choice
of $D\smod$ in the class of $\QQ$-linear equivalence
\textup(respectively $I$-linear equivalence under \eqref{conj-main-adj-21}\textup).
Moreover, if $(X,D)$ is lc and any lc centre of $(X,D)$ dominates $Z$,
then $(Z,D\sdiv+D\smod)$ is klt.
\end{enumerate}
\end{corollary}

\begin{proof}
For a log resolution $g\colon Z'\to Z$ of the pair 
$(Z, D\sdiv)$, consider 
base change \eqref{eq-Diff-sq-diag}.
Thus $\Supp D'\sdiv$ is a simple normal crossing divisor on $Z'$.

(i) 
Put $D_t:=(1-t)D +t \Theta$ (see \xref{assumpt-adj-*}).
Assume that $(X,D)$ is not lc. 
Then $(X,D_t)$ is also not lc for some $0<t\ll 1$.
Let $F$ be a divisor of discrepancy
$a(F,X,D_t)<-1$. Since $(X,D_t)$ is klt near the generic fibre,
the centre of $F$ on $Z$ is a proper subvariety.
By Theorem \ref{th-Ambro-m} we can take $g$
so that 
$\bK+(\mathbf D_t)\sdiv=\ov{K_{Z'}+(D'_t)\sdiv}$\  and 
$(\mathbf D_t)\smod=\ov{(D'_t)\smod}$.
Moreover, by
\cite[Ch. VI, Th. 1.3]{Kollar-1996-RC}
we can also take $g$ so that the centre of $F$ on $Z'$ is a prime divisor,
say $W$. Put $(D_t)_Z:=(D_t)\sdiv+(D_t)\smod$.
By \eqref{eq-Diff-sq} we have
\[
-1\le a(W,Z,(D_t)_Z)=a(W,Z', (D'_t)\sdiv+(D'_t)\smod)\le a(W,Z',(D'_t)\sdiv).
\]
Therefore $(X',D'_t)$ is lc over the generic point of
$W$ (see \eqref{cW}). In particular,
$a(F,X',D'_t)=a(F,X,D_t)\ge -1$, a contradiction.

(ii) 
By our assumption \eqref{conj-main-adj-1} 
$\mathbf D\smod$ is b-Cartier, so
we can take $g$ so that $\mathbf D\smod=\ov{D'\smod}$
and 
$\bK+\mathbf D\sdiv=\ov{K_{Z'}+D'\sdiv}$.
Moreover by \eqref{conj-main-adj-1}
(respectively by \eqref{conj-main-adj-21}) we can take $g$ so that 
${D'\smod}$ (respectively ${ID'\smod}$) 
is semiample (respectively linearly free).
By \eqref{eq-Diff-sq} $g$ is $(K+D\sdiv+D\smod)$-crepant.
If $(Z',D'\sdiv)$ is lc (resp. klt), then 
replacing ${D'\smod}$ with an effective general
representative of the corresponding 
class of $\QQ$-linear equivalence we obtain
\[
\discr(Z,D\sdiv+D\smod)=\discr(Z',D'\sdiv+D'\smod)=\discr(Z',D'\sdiv)
\ge -1.
\]
(resp. $>-1$). 
We can suppose also that $\down{D\smod}=0$.
Hence $(Z,D\sdiv+D\smod)$ is lc (resp. klt) in this case.
Thus we assume that $(Z',D'\sdiv)$ is not lc (resp. not klt).
Let $E$ be a divisor over $Z$ of discrepancy
$a(E,Z',D'\sdiv)\le -1$. Clearly, we may assume that 
$\Center_{Z'} E \not\subset \Supp D'\smod$.
Then $a(E,Z',D'\sdiv+D'\smod)=a(E,Z',D'\sdiv)\le -1$.
Replacing $Z'$ with its blowup we may assume that 
$E$ is a prime divisor on $Z'$ (and again
$\Center_{Z'} E \not\subset \Supp D'\smod$).
Since $(X',D')$ is lc and by \eqref{cW},
$c_E=0$, $d_E=1$, and $a(E,Z',D'\sdiv)=-1$. 
Then $(Z',D'\sdiv)$ is lc.
Furthermore, by \eqref{cW} the pair $(X', D'+cf^{\prime \bullet}E)$
is not lc for any $c>0$. 
This means that $f^{-1}(\Center_Z (E))$ contains an lc centre.
\end{proof}

The following example 
shows that the condition $\mathbf D\smod\ge 0$ in (i) of 
Corollary \ref{cor_Inv_aDj}
cannot be omitted.
\begin{example}
\label{ex-Inv-Adj-contr}
Let $f\colon X\to Z=\CC^2$ be a standard conic bundle
given by $x^2+uy^2+vz^2$ in $\PP^2_{x,y,z}\times \CC^2_{u,v}$.
The linear system 
$ |-n K_X|$ is base point free for $n\ge 1$. 
Let $H\in |-2 K_X|$ be
a general member. 
Now let $\Gamma(t):= \Gamma_1+t\Gamma_2$, where
$\Gamma_1:=\{u=0\}$ and $\Gamma_2:=\{v=0\}$.
Put $D(t):=\frac12 H+f^*\Gamma(t)$. 
Then $2(K+D)= 2f^*\Gamma(t)$ and $D(t)\sdiv=\Gamma(t)$. 
Since $K_Z=0$, we have $D\smod =0$.

For $t=1$,
the log divisor $K_Z+\Gamma(t)$ is lc but
$K+D(t)= f^*(K_Z+\Gamma(t))$ is not. Indeed, in the chart $z\neq
0$ there is an isomorphism 
\begin{equation}
\label{ex-m-cp}
(X,f^*\Gamma)\simeq (\CC^3_{x,y,u}, \{u(x^2+uy^2)=0\}). 
\end{equation}

The explanation of this fact is that the b-divisor
$\bD\smod$ is non-trivial. To show this 
we consider the 
following diagram \cite[\S 2]{Sarkisov-1980-re}:
\[
\xymatrix{
X'\ar@{-->}[r]^{\chi}\ar[dd]^{f'}&\tilde X\ar[dr]^{h}&
\\
&&X\ar[d]^{f}
\\
Z'\ar[rr]^{g}&&Z
}
\]
where $h$ is the blowup the central fibre $f^{-1}(0)_{\red}$,
$\chi$ is the simplest flop, $g$ is the blowup 
of $0$, and $f'$ is again a standard conic bundle.
Put $t=1/2$ and let $\tilde D$ and $D'$ be the crepant 
pull-backs of $D:=D(t)$ on $\tilde X$ and $X'$, respectively.
The $h$-exceptional divisor $F$ appears in $\tilde D$ 
with multiplicity $1/2$. Let $F'$ be the proper transform of $F$ on $X'$.
Then $F'=f'^* E$, where $E$ is the $g$-exceptional 
divisor.
It is easy to see from \eqref{ex-m-cp} that the pair
$(X,D)$ is lc but not klt at the generic point 
of $f^{-1}(0)_{\red}$. So is
$(\tilde X, \tilde D)$ at the generic point
of the flopping curve. This implies 
that $(X', D')$ 
is lc but not klt over the generic point of $E$.
Therefore, 
$D'\sdiv=E+\Gamma'$, where $\Gamma'$ is the proper transform of $\Gamma$.
On $Z'$, we have $K_{Z'}=E$ and $K+D'=f'^*g^*\Gamma$, 
so $D'\smod=g^*\Gamma-E-D'\sdiv=-\frac12 E$.
Thus $D'\smod\le 0$ and $2D'\smod$ is free.
\end{example}

\section{Two important particular cases of Effective Adjunction}
\label{section-Adj-conj-part-cases}
Using a construction and a result of \cite{Kawamata-1997-Adj} we prove the following.
\begin{theorem}
\label{th-n-n-1}
Conjectures \xref{conj-main-adj} hold if $\dim X=\dim Z+1$.
\end{theorem}
\begin{remark}
We expect that in this case one can take
$I=12q$, where $q$ is a positive integer such that 
$qD^{\h}$ is an integral
divisor. 
\end{remark}
\begin{proof}[Proof of Theorem \xref{th-n-n-1}]
We may assume that a general fibre of $f$ is a rational 
curve (see Example \ref{ex-elliptic}).
Thus the horizontal part $D^\h$ of $D$ is non-trivial.
First we reduce the problem to the case when 
all components of $D^\h$ are 
generically sections.
Write $D=\sum d_iD_i$ and take 
\[
\delta:= \min \{ d_i \mid \text{$D_i$ is horizontal and $d_i>0$}\}.
\]
(we allow components with $d_i=0$).
Let $D_i$ be a horizontal component and let 
$D_i\to \hat Z\stackrel{g}{\longrightarrow} Z$ be the Stein factorization of 
the restriction $f|_{D_i}$. Let $n_i:=\deg g$. 
Let $l$ be a general fibre of $f$. 
Since $d_iD_i\cdot l\le D\cdot l=-K\cdot l=2$, we have 
\begin{equation}
\label{eq-Adj-3-2-n-i}
n_i= D_i\cdot l\le 2/d_i\le 2/\delta.
\end{equation}
Assume that $n_i>1$. Consider the base change
\[
\begin{CD}
\hat X@>{h}>>X
\\
@V{\hat f}VV@V{f}VV
\\
\hat Z@>{g}>>Z
\end{CD}
\]
where $\hat X$ is the normalization of the dominant component of $X\times_{Z} \hat Z$.
Define $\hat D$ on $\hat X$ by 
\begin{equation}
\label{eq-3-1-crepant-pb}
K_{\hat X}+\hat D=h^*(K_X+D).
\end{equation}
More precisely, $\hat D=\sum_{i,j} \hat d_{i,j}\hat D_{i,j}$, where
$h(\hat D_{i,j})=D_i$, $1-\hat d_{i,j}=r_{i,j}(1-d_i)$, and $r_{i,j}$
is the ramification index along $\hat D_{i,j}$. By construction,
the ramification locus $\Lambda$ of $h$ is $\hat f$-exceptional, that is 
$\hat f (\Lambda)\neq \hat Z$.
Therefore, $\hat D$ is a boundary near the generic fibre.
Similarly, we define $\hat \Theta$ as the crepant pull-back of $\Theta$ 
from \ref{assumpt-adj-*}. Thus the pair $(\hat X,\hat D)$ satisfies
assumptions of \xref{not-adj} and \xref{assumpt-adj-*}.
It follows from \eqref{eq-3-1-crepant-pb} that 
\[
K_{\hat X}+\hat D=\hat f^* g^*(K_Z+D\sdiv+D\smod).
\]
According to \cite[Th. 3.2]{Ambro-PhD} for the discriminant $\hat D\sdiv$ of 
$\hat f$ we have
\[
K_{\hat Z}+\hat D\sdiv=g^*(K_Z+D\sdiv).
\]
For a suitable choice of $\hat D\smod$ in the class 
of $n_i$-linear equivalence, we can write
\[
K_{\hat Z}+\hat D\sdiv+\hat D\smod=g^*(K_Z+D\sdiv+D\smod).
\]
Therefore, $\hat D\smod=g^* D\smod$. 
If $\hat I \hat D\smod$ is free for some positive integer $\hat I$, then
so is $n_i \hat I D\smod$.
Thus we have proved the following.
\begin{claim}
\label{claim-Adj-3-2-a}
Assume that Conjecture \xref{conj-main-adj} holds for 
$\hat f\colon \hat X\to \hat Z$ with constant $\hat I$.
Then this conjecture holds for 
$f\colon X\to Z$ with $I:=n_i \hat I$.
\end{claim}

Note that the restriction
$\hat f|_{h^{-1} (D_i)}\colon h^{-1} (D_i)\to \hat Z$ 
is generically finite of degree $n_i$. Moreover, 
$h^{-1} (D_i)$ has a component which is a section
over the generic point. 
Applying Claim \ref{claim-Adj-3-2-a} several times and taking 
\eqref{eq-Adj-3-2-n-i} into account we obtain the desired reduction 
to the case when all the horizontal $D_i$'s 
with $d_i>0$ are generically sections. 

\subsection{}
\label{subs-Adj-conic-sect}
Further by making a birational base change and by blowing up $X$
we can get the situation when 
\begin{enumerate}
\item 
$Z$ and $X$ are smooth, 
\item 
the $D_i$'s are regular disjointed sections,
\item 
the morphism 
$f$ is smooth outside of a simple normal crossing divisor
$\Xi\subset Z$,
\item 
$f^{-1}(\Xi)\cup \Supp D$ is also a simple normal crossing divisor.
\end{enumerate}
Let $n$ be the number of horizontal components of $D$.
Note that we allow sections with multiplicities $d_i=0$ on this step.

Let $\mathcal M_n$ be the moduli space of 
$n$-pointed stable rational curves, let 
$f_n\colon \mathcal U_n\to \mathcal M_n$ 
be the corresponding universal family, and let 
$\mathcal P_1,\dots,\mathcal P_n$ be sections of $f_n$ which correspond to 
the marked points (see \cite{Knudsen-1983-2}).
It is known that both $\mathcal M_n$ and $\mathcal U_n$ are smooth and projective.
Take $d_i\in [0,\, 1]$ so that $\sum d_i=2$ and 
put $\mathcal D:=\sum d_i\mathcal P_i$. Then $K_{\mathcal U_n}+\mathcal D$
is trivial on the general fibre. 
However, $K_{\mathcal U_n}+\mathcal D$ is not 
numerically trivial everywhere over $\mathcal M_n$, moreover, it is not nef everywhere 
over $\mathcal M_n$: 

\begin{theorem}[see \cite{Keel-1992}, \cite{Kawamata-1997-Adj}]
\label{th-Knud}
\begin{enumerate}
\item 
There exist a smooth projective variety $\bar {\mathcal{U}}_n$, a 
$\PP^1$-bundle $\bar f_n\colon \bar {\mathcal{U}}_n\to \mathcal M_n$, and a sequence of 
blowups \textup(blowdowns\textup) with smooth centres 
\[
\sigma\colon {\mathcal{U}}_n= \mathcal{U}^{1}\to \mathcal{U}^{2}\to \cdots\to 
\mathcal{U}^{n-2} = \bar {\mathcal{U}}_n.
\]
\item 
For $\bar{\mathcal D}:=\sigma_*\mathcal D$, the \textup(discrepancy\textup) 
divisor 
\[
\mathcal F:=K_{{\mathcal{U}}_n}+\mathcal D-\sigma^*(K_{\bar {\mathcal{U}}_n}+\bar {\mathcal D})
\]
is effective and \emph{essentially exceptional} 
on $\mathcal M_n$. 
\item
There exists a semiample $\QQ$-divisor $\mathcal L$ on $\mathcal M_n$ such that
\[
K_{\bar {\mathcal{U}}_n}+\bar {\mathcal D}=
\bar f_n^*(K_{\mathcal M_n}+\mathcal L).
\] 
Therefore,
\[
K_{{\mathcal{U}}_n}+\mathcal D -\mathcal F= 
f_n^*(K_{\mathcal M_n}+\mathcal L).
\]
\end{enumerate}
\end{theorem}
Recall that
for any contraction $\var\colon Y\to Y'$,
a divisor $G$ on $Y$ is said to be 
\textit{essentially exceptional} over $Y'$ if 
for any prime divisor $P$ on $Y'$, the support of the divisorial pull-back
$\var^\bullet P$ is not contained in $\Supp G$.

\begin{corollary}
In the above notation we have 
\[
(\mathcal D-\mathcal F)\sdiv=0,\qquad 
(\mathcal D-\mathcal F)\smod=\mathcal L.
\]
\end{corollary}
Moreover,
the proof of Theorem \ref{th-n-n-1} implies that the b-divisor
$\mathbf G$, the b-divisor of the moduli part of $\mathcal D-\mathcal F$,
stabilizes on $\mathcal M_n$, that is,
$\mathbf G=\ov{(\mathcal D-\mathcal F)\smod}$. 
\begin{proof}
See Example \ref{exam-mod-mod} below.
\end{proof}

Since the horizontal components of $D$ 
are sections, $(X/Z,D^{\h})$
is generically an $n$-pointed stable curve
\cite{Knudsen-1983-2}.
Hence we have the induced rational maps
\[
\xymatrix{
X\ar@{-->}[r]^{\beta}\ar[d]^{f}&\mathcal U_n\ar[d]^{f_n}
\\
Z\ar@{-->}[r]^{\phi}&\mathcal M_n
}
\]
so that $f_n\comp \beta=\phi \comp f$ 
and $\beta (D_i)\subset \mathcal P_i$.
Let $\Xi\subset Z$ as above (see \ref{subs-Adj-conic-sect}, (iii)).
Thus $f$ is a smooth morphism over $Z\setminus \Xi$.
Replacing $X$ and $Z$ with its birational models
and $D$ with its crepant pull-back
we may assume additionally to \ref{subs-Adj-conic-sect}
that $\beta$ and $\phi$ are
regular morphisms.
Now take the $\mathcal D=\sum d_i\mathcal P_i$ so that 
it corresponds to the horizontal part 
$D^{\h}=\sum_{f(D_i)=Z} d_i D_i$.
Consider the following commutative diagram
\[
\xymatrix{
X\ar[r]^{\mu}\ar@/^2pc/[rrr]^{\beta}
\ar[ddr]^{f}
&\hat X\ar[rr]^{\psi}\ar[dd]^{\hat f}&&{\mathcal U}_n
\ar[dd]^{f_n}\ar@/^0.3pc/[dr]^{\sigma}
\\
&&&&\bar {{\mathcal U}_n}\ar@/^0.3pc/[dl]^{\bar f_n}
\\
&Z\ar[rr]^{\phi}&&{\mathcal M}_n
}
\]
where $\hat X:= Z\times_{\mathcal M_n} {\mathcal U}_n$.

\begin{edefinition}
\label{Adj-conic-hat-plt}
Since the fibres of $f_n$ are stable curves, 
near every point $u\in \mathcal U_n$
the morphism $f_n$ is either smooth or in a suitable local
analytic coordinates is given by 
\[
(u_1,\, u_2,\dots,u_{n-2}) \longmapsto 
(u_1u_2,\, u_2,\dots,u_{n-2}).
\] 
Then easy local computations show 
that $\hat X$ is normal and has only canonical singularities
\cite{Kawamata-1997-Adj}.
Moreover, the pair $(\hat X,\hat D^{\h}=\psi^*\mathcal D)$
is canonical because $f_n$ is a smooth morphism near $\Supp \mathcal D$. 
\end{edefinition}
We have 
\[
K_{X}+D= f^*(K_{Z}+D\sdiv+D\smod).
\]
Put $\hat D:=\mu_*D$. Then $K_{X}+D=\mu^*(K_{\hat X}+\hat D)$, 
so 
\[
\hat D\sdiv=D\sdiv, \qquad \hat D\smod=D\smod,
\]
\[
K_{\hat X}+\hat D= \hat f^*(K_{Z}+D\sdiv+D\smod).
\]

\begin{edefinition}
Let $\var\colon Y\to Y'$ be any contraction, where 
$\dim Y'\ge 1$.
We introduce $G^{\bot}= G- \var^\bullet G_{\neg}$, where 
$G_{\neg}$ is taken so that
the vertical part $(G^{\bot})^\ver$ of $G^{\bot}$ is
essentially exceptional and $G_{\neg}$ is maximal with this property.
In particular, $(G^{\bot})^\ver\le 0$ over an open subset $U'\subset Y'$
such that $\codim (Y'\setminus U')\ge 2$. 
Note that our construction of $G_{\neg}$ and 
$G^{\bot}$ is in codimension one over $Y'$, i.e., 
to find $G_{\neg}$ and 
$G^{\bot}$ we may replace $Y'$ with $Y'\setminus W$,
where $W$ is a closed subset of codimension $\ge 2$.
\end{edefinition}

\begin{lemma}
\label{lemma-adj-minimal}
Let $\var\colon Y\to Y'$ be a contraction 
and let $G$ be an $\RR$-divisor on $Y$.
Assume that $\dim Y'\ge 1$.
Assume that $(Y/Y',G)$ satisfies conditions \xref{not-adj}.
The following are equivalent:
\begin{enumerate}
\item 
$G^{\ver}-\var^\bullet G\sdiv$ 
is essentially exceptional,
\item 
$G\sdiv=G_\neg$,
\item 
$(G^\bot)\sdiv=0$.
\end{enumerate}
\end{lemma}

\begin{proof}
Implications (ii) $\Longleftrightarrow$ (iii) $\Longrightarrow$ (i)
follows by definition of $G\sdiv$ and 
semiadditivity (Lemma \ref{lemma-adj-prelim}).
Let us prove (i) $\Longrightarrow$ (ii).
Assume that $G^{\ver}-\var^\bullet G\sdiv$ 
is essentially exceptional.
Then by definition $G\sdiv\le G_\neg$.
On the other hand, for any prime divisor $P\subset Y'$, the multiplicity 
of $G^{\bot}$ along some component of $\var^\bullet P$
is equal to $0$.
Hence the log
canonical threshold of $(K+G^\bot, \var^\bullet P)$ 
over the generic point of $P$ is $\le 1$.
So by definition of the 
divisorial part and Lemma \ref{lemma-adj-prelim} we have 
$0\le (G^{\bot})\sdiv=G\sdiv-G_\neg$.
\end{proof}

\begin{example}
\label{exam-mod-mod}
Clearly, for $f_n\colon \mathcal U_n \to \mathcal M_n$,
the discrepancy divisor
$\mathcal F$ is essentially exceptional.
Hence,
$(\mathcal D-\mathcal F)\sdiv\ge 0$.
On the other hand, 
by construction every fibre of $f_n$ is reduced.
Hence, for every prime divisor $W\subset Z$,
the divisorial pull-back $f_n^\bullet W$ is reduced and 
$\Supp (f_n^\bullet W + \mathcal D)$ is a simple normal 
crossing divisor over the generic point of $W$.
This implies that $c_W\ge 1$ and so $(\mathcal D-\mathcal F)\sdiv=0$.
\end{example}

\begin{proof}[Proof of Theorem {\xref{th-n-n-1}} \textup(continued\textup)]
It is sufficient to show that 
$D\smod=\phi^*\mathcal L=\phi^*(\mathcal D-\mathcal F)\smod$
(we replace $Z$ with its blowup if necessary).
Then the b-divisor $\mathbf D\smod$ automatically
stabilizes on $Z$, i.e., $\mathbf D\smod=\ov{D\smod}$.
In this situation $D\smod$ is effectively semiample 
because $N\mathcal L$ is an integral base point free 
divisor for some $N$ which depends only on $n$.
Since and $\phi$ is a regular morphism,
to show $D\smod=\phi^*\mathcal L$ we will
freely replace $Z$ with an open subset $U\subset Z$
such that $\codim (Z\setminus U)\ge 2$.
Thus all the statements below are valid
over codimension one over $Z$. In particular,
we may assume that $D\smod=(D^\bot)\smod$.
Replacing $D$ with $D^\bot$ we may assume that 
$D_\neg=0$ (we replace $Z$ with $U$ as above).
Thus $D^\ver\le 0$ 
and $D^\ver$ is essentially exceptional.
In particular, $D\sdiv \ge 0$. 

On the other hand, by construction the fibres 
$(\hat f^*(z),\hat D^{\h}=\psi^*\mathcal D)$, $z\in Z$
are stable (reduced) curves. In particular, they are 
slc (semi log canonical \cite[\S 4]{Kollar-ShB-1988},
\cite[Ch. 12]{Utah}). 
By the inversion of adjunction \cite[\S 3]{Shokurov-1992-e},
\cite[Ch. 16-17]{Utah}
for every prime divisor $W\subset Z$ 
and generic hyperplane sections $H_1,\dots,H_{\dim Z-1}$
the pair 
$(\hat X, \hat D^{\h}+
\hat f^\bullet W+\hat f^\bullet H_1+\cdots+\hat f^\bullet H_{\dim Z-1})$
is lc. Since $\hat D^{\ver}\le 0$, so is 
the pair $(\hat X, \hat D+\hat f^\bullet W)$.
This implies that $c_W\ge 1$ and so $D\sdiv=\hat D\sdiv=0$.

We claim that $\hat D^{\ver}=\mu_*D^{\ver}$ 
is essentially exceptional. Indeed, otherwise 
$\hat D^{\ver}$ is strictly negative over
the generic point of some prime divisor $W\subset Z$,
i.e., $\mu$ contracts all the components $E_i$ of $f^\bullet W$ 
of multiplicity $0$. By \ref{Adj-conic-hat-plt} the pair 
$(\hat X,\hat D+\ep \hat f^\bullet W)$ is canonical over the generic 
point of $W$ for some small positive $\ep$.
On the other hand, for the
discrepancy of $E_i$ we have
$a(E_i,\hat X,\hat D+\ep \hat f^\bullet W)=
a(E_i,X,D+\ep f^\bullet W)=-\ep$. The contradiction proves our 
claim.

For relative canonical divisors we have 
\[
K_{\hat X/Z}=\psi^*K_{{\mathcal U}_n/{\mathcal M}_n}
\]
(see, e.g., \cite[Ch. II, Prop. 8.10]{Hartshorn-1977-ag}).
Taking
$\hat D^{\h}=\psi^*\mathcal D$ into account we obtain
\[
K_{\hat X/Z}+\hat D^{\h}- \psi^*\mathcal F=\psi^*
(K_{{\mathcal U}_n/{\mathcal M}_n}+\mathcal D-\mathcal F)
=\psi^* f_n^* \mathcal L=\hat f^* \phi^* \mathcal L.
\]
Hence,
\[
-\hat D^{\ver}-\psi^*\mathcal F\q{\RR}
K_{\hat X}+\hat D^{\h}-\psi^*\mathcal F
\q{\RR} \hat f^* \phi^* \mathcal L+\hat f^*K_Z
\]
over $Z$, i.e., $\hat D^{\ver}+\psi^*\mathcal F$ is 
$\RR$-linearly
trivial over $Z$.

Since $\psi^*\mathcal F$ is also 
essentially exceptional over 
$Z$, by Lemma \ref{11a} below we have 
$\hat D^{\ver}= - \psi^*\mathcal F$ and
\[
\hat f^*D\smod=\hat f^*(D\smod+D\sdiv)=K_{\hat X/Z}+\hat D=\hat f^* \phi^* \mathcal L.
\]
This gives us $D\smod=\phi^*\mathcal L=\phi^*(\mathcal D-\mathcal F)\smod$.
Therefore $D\smod$ is effectively semiample.
This proves Theorem \ref{th-n-n-1}.
\end{proof}

\begin{lemma}[cf. {\cite[Lemma 1.6]{Prokhorov-2003d-e}}]
\label{11a}
Let $\var\colon Y\to Y'$ be a contraction with $\dim Y'\ge 1$
and let $A$, $B$ be essentially exceptional over $Y'$ 
divisors on $Y$ such that 
$A\equiv B$ over $Y'$ and $A,\, B\le 0$
\textup(both conditions are over codimension one over $Y'$\textup).
Then $A=B$ over codimension one over $Y'$.
\end{lemma}

\begin{proof}
The statement is well-known in the birational case (see
\cite[\S 1.1]{Shokurov-1992-e}), so we assume that 
$\dim Y'<\dim Y$. As in \cite[Lemma 1.6]{Prokhorov-2003d-e},
replacing $Y'$ with its general hyperplane
section $H'\subset Y'$ and $Y$ with $\var^{-1} (H')$ we 
may assume that $\dim \var (\Supp A) =0$
and $\dim \var (\Supp B) \ge 0$.
The essential exceptionality of $A$ and $B$ is preserved. 

We may also assume that
$Y'$ is a sufficiently small affine neighbourhood of some fixed
point $o\in Y'$ (and $\var (\Supp A) =o$). Further, all the
conditions of lemma are preserved if we replace $Y$ with its
general hyperplane section $H$. If $\dim Y' > 1$, then we can
reduce our situation to the case $\dim Y=\dim Y'$. Then the
statement of the lemma follows by
\cite[\S 1.1]{Shokurov-1992-e}
and from the
existence of the Stein factorization. 
Finally, consider the case
$\dim Y'=1$ (here we may assume that $\dim Y=2$
and $\var$ has connected fibres). 
By the Zariski
lemma $A=B+a\var^*o$ for some $a\in \QQ$.
Since $A$ and $B$ are essentially exceptional
and $\le 0$, $a=0$.
\end{proof}

\end{proof}

\begin{example}
Assume that all the components $D_1,\dots,D_r$ 
of $D^{\h}$ are sections.
If $r=3$, then since $\mathcal M_3$ is a point, we have
$D\smod=0$. For $r=4$ the situation is more complicated:
$\mathcal M_4\simeq\PP^1$,\quad $\mathcal U_4$ is a del Pezzo surface 
of degree $5$, and $f_4\colon \mathcal U_4\to \mathcal M_4=\PP^1$
is a conic bundle with three degenerate fibres. Each component
of degenerate fibre meets exactly two components of $\mathcal D$.
Hence $\bar {\mathcal D}$ is a normal crossing divisor. 
It is easy to see that $\sigma$ contracts a
component of a degenerate fibre which meets 
$\mathcal D_i$ and $\mathcal D_j$ with $d_i+d_j\le 1$.
Clearly, $\bar {\mathcal U}_4\simeq \mathbb F_e$ 
is a rational ruled surface, $e= 0$ or $1$.
We can write $\bar {\mathcal D}_i\sim \Sigma+a_i F$,
where $\Sigma$ is the minimal section and 
$F$ is a fibre of $\bar {\mathcal U}_4=\mathbb F_e\to \PP^1$.
Up to permutation we may assume that $\bar {\mathcal D}_i\neq \Sigma$
for $i=2,3,4$. Taking $\sum d_i=2$ into account we get
\[
K_{\bar {\mathcal U}_4}+\bar {\mathcal D}
\sim -2\Sigma-(2+e)F+\sum d_i(\Sigma+a_i F)=
\left( \sum d_i a_i-e\right) F+\bar f_n^*K_{\mathcal M_4}. 
\]
Therefore, 
\[
\deg \mathcal L= \sum d_i a_i-e\ge e\sum d_i- ed_1-e\ge 0.
\]
\end{example}

\subsection{}
Now we consider the case when the base variety $Z$ is a curve.

\begin{proposition}
\label{Prop-dimz=1}
Assume Conjectures~\xref{BAB} and \xref{conj-main-adj}
in dimensions $\le d-1$ and
LMMP in dimension $\le d$. 
If $X$ is FT \textup(and projective\textup) variety of dimension $d$, 
then Conjecture~ \xref{conj-main-adj} holds in dimension $d$.
\end{proposition}

\begin{corollary}
\label{corollary-815}
Conjecture \xref{conj-main-adj} holds true in the following cases\textup:
\begin{enumerate}
\item
$\dim X=\dim Z+1$,
\item
$\dim X=3$ and $X$ is FT.
\end{enumerate}
\end{corollary}
\begin{proof}
Immediate by Theorem \ref{th-n-n-1}
and Proposition \ref{Prop-dimz=1}.
\end{proof}

The rest of this section is devoted to the proof of 
Proposition \ref{Prop-dimz=1}.
Thus from now on and through the end of this section 
we assume that the base variety $Z$ is a curve.
First we note that $Z\simeq \PP^1$ because 
$X$ is FT.

\begin{lemma}
\label{lemma-local-comlpl-31}
Fix a positive integer $N$.
Let $f\colon X\to Z\ni o$ be a contraction to a curve germ
and let $D$ be an 
$\RR$-divisor on $X$. 
Let $D^{\h}$ be the horizontal part of $D$. Assume that
\begin{enumerate}
\item 
$\dim X\le d$ and $X$ is FT over $Z$,
\item
$D^\h$ is a $\QQ$-boundary and $ND^{\h}$ is integral, 
\item
$K_X+D$ is lc and numerically trivial over $Z$. 
\end{enumerate}
Assume LMMP in dimension $\le d$.
Further assume that the statement of Theorem
\xref{main-result} holds in dimensions $\le d-1$.
Then there is an $n$-complement $K+D^+$ 
of $K+D$ near $f^{-1}(o)$
such that $N \mid n$, $n\le \Const (N,\dim X)$,
and $a(E,X,D^+)=-1$ for some divisor $E$ with 
$\Center_Z E=o$.
\end{lemma}
\begin{proof}
Take a finite set $\R\subset [0,\, 1]\cap \QQ$ 
and a positive integer $I$ so that $D^{\h}\in \R$, $I(\R)\mid I$, and $N\mid I$. 
Replacing $D$ with $D+\alpha f^*o$ we may assume that 
$(X,D)$ is maximally lc. 
Next replacing $(X,D)$ with its suitable blowup
we may assume that $X$ is 
$\QQ$-factorial and 
the fibre $f^{-1}(o)$ has 
a component, say $F$, of multiplicity $1$
in $D$.
Run $-F$-MMP over $Z$. 
This preserves the $\QQ$-factoriality and lc property of 
$K+D$. Clearly, $F$ is not contracted. 
On each step, the contraction is birational.
So at the end we get a model with irreducible central fibre: $f^{-1}(o)_{\red}=F$.
Then $D\in \Phi(\R)$.
Applying $D^{\h}$-MMP over $Z$,
we may assume that $D^{\h}$ is nef over $Z$.
We will show that $K+D$ is $n$-complemented for some $n\in \NNN_{d-1}(\ov \R)$.
Then by Proposition \ref{cor-pull-back_compl-I} we can pull-back complements 
to our original $X$. 
Note that the $f$-vertical part of $D$ coincides with $F$,
so it is numerically trivial over $Z$.
Since $X$ is FT over $Z$, $-K_X$ is big over $Z$. 
Therefore $D\equiv D^{\h}$ is nef and big over $Z$. 
Now apply construction of \cite[\S 3]{Prokhorov-Shokurov-2001} to
$(X, D)$ over $Z$.
There are two cases:
\begin{enumerate}
\item[(I)] 
$(X,F)$ is plt,
\item[(II)]
$(X,F)$ is lc but not plt (recall that $F\le D$).
\end{enumerate}
Consider, for example, the second case (the first case
is much easier and can be treated in a similar way).
First we define an auxiliary boundary to localise a 
suitable divisor of discrepancy $-1$. 
By Kodaira's lemma, for some effective $D^\mho$, 
the divisor $D-D^\mho$ is ample. 
Put $D_{\ep,\alpha}:=(1-\ep)D+\alpha D^\mho$.
Then $K_X+D_{\ep,\alpha}\equiv -\ep D+\alpha D^\mho$.
So $(X,D_{\ep,\alpha})$ is a klt log Fano over $Z$
for $0<\alpha \ll \ep \ll 1$.
Take $\beta=\beta(\ep,\alpha)$ so that 
$(X,D_{\ep,\alpha}+\beta F)$ is maximally lc
and put $G_{\ep,\alpha}:=D_{\ep,\alpha}+\beta F$.
Thus $(X,G_{\ep,\alpha})$ is a lc (but not klt) log Fano over $Z$.

Let $g\colon \widehat X\to X$ be an \textit{inductive blowup} of $(X,G_{\ep,\alpha})$
\cite[Proposition 3.6]{Prokhorov-Shokurov-2001}.
By definition $\widehat X$ is $\QQ$-factorial, 
$\rho(\widehat X/X)=1$, the $g$-exceptional locus is a prime divisor
$E$ of discrepancy $a(E,X,G_{\ep,\alpha})=-1$, the pair
$(\widehat X,E)$ is plt, and $-(K_{\widehat X}+E)$
is ample over $X$. 
Since $(X,G_{\ep,\alpha}-\gamma F)$ is klt for $\gamma>0$, 
$\Center_Z(E)=o$.
Note that, by construction, $E$ is not
exceptional on some fixed log resolution of $(X,\Supp G_{\ep,\alpha})$.
Hence we may assume that $E$ and $g$ do not depend on 
$\ep$ and $\alpha$ if $0<\ep \ll 1$.
In particular, $a(E,X,D)=-1$.

By (iii) of Lemma \ref{lemma-FT} $\widehat X$ is FT over $Z$. 
Let $\widehat D$ and $\widehat G_{\ep,\alpha}$
be proper transforms on $\widehat X$ of $D$ and $G_{\ep,\alpha}$, respectively.
Then
\[
\begin{array}{rll}
0\equiv &g^*(K_X+D)&=K_{\widehat X}+\widehat D+E,
\\[8pt]
 &g^*(K_X+G_{\ep,\alpha})&=K_{\widehat X}+\widehat G_{\ep,\alpha}+E,
\end{array}
\]
where $-(K_X+G_{\ep,\alpha})$ is ample over $Z$.
Run $-(K_{\widehat X}+E)$-MMP starting from $\widehat X$ over $Z$:
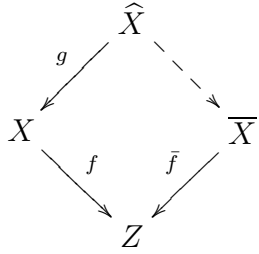
\[
\xymatrix{
&\widehat X\ar[dl]_{g}\ar@{-->}[dr]^{}&
\\
X\ar[dr]^{f}&&\ov X\ar[dl]_{\bar f}
\\
&Z&
}
\]
Since $-(K_{\widehat X}+E)\equiv \widehat D$, we can contract only components 
of $\widehat D$. 
At the end we get a model $(\ov X,\ov D+\ov E)$ such that
$-(K_{\ov X}+\ov E)$ is nef and big
over $Z$, $K_{\ov X}+\ov E+\ov D\equiv 0$, 
and $(\ov X,\ov E+\ov D)$ is lc.

We claim that the plt property of $K_{\widehat X}+E$ is preserved 
under this LMMP.
Indeed, for $0<t\ll 1$, the log divisor 
$K_{\widehat X}+(1-t)\widehat G_{\ep,\alpha}+E$ 
is a convex linear combination of log divisors
$K_{\widehat X}+\widehat G_{\ep,\alpha}+E$ 
and
$K_{\widehat X}+E$.
The first divisor is anti-nef and is trivial only on one 
extremal ray $R$, the ray generated by fibres of $g$.
The second one is strictly negative on $R$.
Since $\widehat X$ is FT over $Z$, the Mori cone $\NE(\widehat X/Z)$
is polyhedral.
Therefore $K_{\widehat X}+(1-t)\widehat G_{\ep,\alpha}+E$
is anti-ample (and plt) for $0<t\ll 1$.
By the base point free theorem there is a boundary 
$M\ge (1-t)\widehat G_{\ep,\alpha}+E$ such that 
$(\widehat X,M)$ is a plt $0$-pair. Since $E$ is not contracted,
this property is preserved under our LMMP.
Hence $(\ov X,\ov M)$ is plt and so is $(\ov X,\ov E)$.
This proves our claim.
In particular, $\ov E$ is normal and FT.

Take $\delta:=1/m$, $m\in \ZZ$,
$m\gg 0$.
For any such $\delta$, the pair $(\ov X,(1-\delta)\ov D+\ov E)$ is plt
and $-(K_{\ov X}+(1-\delta)\ov D+\ov E)$ is nef and big over $Z$.
By our inductive hypothesis there is an $n$-complement $K_{\ov E}+\Diff_{\ov E}(\ov D)^+$
of $K_{\ov E}+\Diff_{\ov E}(\ov D)$ with $n\in \NNN_{d-1}(\ov \R)$.
Clearly, this is also an $n$-complement of
$K_{\ov E}+\Diff_{\ov E}((1-\delta)\ov D)$.
Note that $nD$ is integral.
We claim that $(1-\delta)\ov D\in \PPP_n$.
Indeed, the vertical multiplicities of $(1-\delta)\ov D$
are contained in $\Phi(\R)$. Let $d_i$ be the
multiplicity of a horizontal component of $\ov D$. Then $nd_i\in \ZZ$.
If $d_i=1$, then obviously $(1-\delta)d_i\in \PPP_n$. 
So we assume that $d_i<1$.
Then $\down{(n+1)d_i}=nd_i$ and
$\down{(n+1)(d_i-\delta)}=nd_i\ge n(d_i-\delta)$ for $\delta\ll 1$.
This proves our claim. Now the same arguments
as in \cite[\S 3]{Prokhorov-Shokurov-2001} shows that
$K_X+(1-\delta)D$ is $n$-complemented near
$f^{-1}(o)$. Since $D\in \PPP_n$, there is
an $n$-complement $K_X+D^+$ of $K_X+D$ near
$f^{-1}(o)$ and moreover, $a(E,X,D^+)=-1$.
\end{proof}

\begin{corollary}
\label{cor-mult-fs}
Notation as in Proposition \xref{Prop-dimz=1}.
The multiplicities of $D\smod$ are contained in a finite set.
\end{corollary}

\begin{proof}
Consider a local $n$-complement $D^+$ of $K+D$ near $f^{-1}(o)$. 
Then $n(K_Z+D^+\sdiv+D^+\smod)$ is integral at $o$.
By construction, $(X,D^+)$ has a centre of log canonical singularities 
contained in $f^{-1}(o)$. Hence $D^+\sdiv=0$. 
By semiadditivity (see Lemma \ref{lemma-adj-prelim}) we have 
$D^+\smod=D\smod$. Thus $nD\smod$ is integral at $o$.
\end{proof}

\begin{proof}[Proof of Proposition \xref{Prop-dimz=1}]
The statement of \eqref{conj-main-adj-1} follows by 
Theorem \ref{th-Ambro-m}
(cf. \cite{Kawamata-1998}). Indeed, for any $0<t<1$ we 
put $D_t:=(1-t)D+t\Theta$, where $\Theta$ is such as in 
\xref{assumpt-adj-*}. Then by 
Theorem \ref{th-Ambro-m}
$(D_t)\smod$ is semiample. Hence so is $D\smod$.

Assertion \eqref{conj-main-adj-21} 
follows by 
Theorem \ref{main-result0} (in lower dimension). 

Finally for \eqref{conj-main-adj-2} we note that 
by Corollary \ref{cor-mult-fs} 
$ID\smod$ is integral and base point free for a bounded $I$
because $Z\simeq \PP^1$.
\end{proof}

\begin{remark}
It is possible that 
Proposition \xref{Prop-dimz=1} can be proved 
by using results of \cite{Fujino-Mori-2000},
\cite{Fujino-2003-AG}. 
In fact, in these papers the authors write down 
the canonical bundle formula (for arbitrary $\dim Z$) 
in the following form
(we change notation a little): 
\[
b(K+D)= f^*(bK_Z+L_{X/Z}^{log, ss})+ \sum_P s_P^D f^*P+B^D.
\]
Here $D\sdiv=\frac1b \sum_P s_P^D P$,
$D\smod=\frac 1b L_{X/Z}^{log, ss}$,
and $\codim f(B^D)\ge 2$, so
the term $B^D$ is zero in our situation.
Under the additional assumption that $D$ is a boundary
it is proved that the denominators of $D\smod$
are bounded (and $D\smod$ is semiample because it is nef
on $Z=\PP^1$), 
see 
\cite[Theorem 4.5]{Fujino-Mori-2000},
\cite[Theorem 5.11]{Fujino-2003-AG}.
This should imply our Proposition \ref{Prop-dimz=1}.
We however do not know 
how to avoid the effectivity condition of $D$.
\end{remark}

\section{The main theorem: Case $-(K+D)$ is nef}
\label{sec-nef}
In this section we prove Theorem \ref{main-result} in case \eqref{case-nef}
and Theorem \ref{main-result0} in the case when $(X,B)$ is not klt.
Thus we apply reduction from \S \ref{sect-reduction} and replace 
$(X,B)$ with $(Y,B_Y)$ and put $D:=D_Y$.
The idea of the proof is to consider the 
contraction $f\colon X\to Z$ 
given by $-(K+D)$ and
use Effective Adjunction to pull-back complements from $Z$. 
In practice, there are several technical issues 
which do not allow us to weaken the last assumption in 
Theorem \ref{main-result},
that is, we cannot omit the klt condition when $K+B\not \equiv 0$. 
Roughly speaking
the inductive step work if the following two 
conditions hold:
\begin{enumerate}
\item 
$0<\dim Z<\dim X$, and
\item
the pair $(Z,D\sdiv+D\smod)$ satisfies 
assumptions of Theorem \ref{main-result}.
\end{enumerate}
The main technical step of the proof is Proposition \ref{prop-D-not-big}.
The proof is given in 
\ref{subsetion-0-pairs-proof} and
\ref{subsect-1-1}.

\subsection{Setup}
\label{not-last-first}
Let $(X,D)$ be an lc log pair
and let $f\colon X\to Z$ be a contraction such that
$K+D \q{\QQ}f^* L$ for some $L$ and
$X$ is FT. 
Further, assume the LMMP in dimension $d:=\dim X$.
Our proof uses induction by $d$.
So we also assume that Theorems \ref{main-result}
and \ref{main-result0}
hold true for all $X$ of dimension $<d$.

By Lemma \ref{lemma-Adj-hor-mult-1} we have the following.

\begin{corollary}
\label{cor-Adj-hor-mult} 
In notation of \xref{not-last-first} assume that $\dim Z>0$.
Fix a finite rational set $\R\subset [0,1]$
and let $D\in \Phi(\R)$.
Then the multiplicities of horizontal components of $D$ 
are contained into a finite subset 
$M\subset \Phi(\R)$, where $M$ depends only on 
$\dim X$ and $\R$.
\end{corollary}

\begin{proof}
Restrict $D$ to a general fibre and apply Lemma 
\ref{lemma-Adj-hor-mult-1}.
\end{proof}

Now we verify that under certain assumptions and conjectures
the hyperstand multiplicities transforms to hyperstandard ones after
adjunction. 

For a subset $\R\subset [0,1]$, denote
\[
\R(n):=\left(\ov \R+\frac 1n\ZZ\right)\cap [0,1],
\qquad \R':= \bigcup_{n\in \NNN_{d-1}(\ov \R)} \R(n)\subset [0,\, 1].
\]
These sets are rational and finite whenever so is $\R$.

\begin{proposition}
\label{prop-R-prime}
In notation of \xref{not-last-first}, 
fix a finite rational set $\R\subset [0,1]$.
\begin{enumerate}
\item 
If $D\in \Phi(\R)$, then $D\sdiv\in \Phi(\R')$.
\item 
If $D\in \Phi(\R, \ep_{d-1})$, then
$D\sdiv\in \Phi(\R',\ep _{d-1})\subset \Phi(\R',\ep _{d-2})$.
\end{enumerate}
\end{proposition}

\begin{proof}
By taking general hyperplane sections we may assume that $Z$ is a
curve. Furthermore, we may assume that $X$ is $\QQ$-factorial.
Fix a point $o\in Z$.
Let $d_o$ be the multiplicity of $o$ in $D\sdiv$. Then $d_o=1-c_o$, where
$c_o$ is computed by \eqref{cW}.
It is sufficient to show that 
$d_o\in \Phi(\R(n)) \cup [1-\ep _{d-1},\, 1]$
for any point $o\in Z$ and some $n\in \NNN_{d-1}(\ov \R)$.
Clearly, we can consider $X$ and $Z$ small neighbourhoods of 
$f^{-1}(o)$ and $o$,
respectively.
We also may assume that $c_o>0$, so $f^{-1}(o)$ does not contain any
centres of log canonical singularities of $(X,D)$.
By our assumptions in \ref{not-last-first}
and Lemma \ref{lemma-local-comlpl-31} there is an $n$-complement 
$K_X+D^+$ of $K_X+D$ near
$f^{-1}(o)$ with $n\in \NNN_{d-1}(\ov \R)$ and moreover, $a(E,X,D^+)=-1$.

Now we show that $d_o\in \Phi(\R(n)) \cup [1-\ep _{d-1},\, 1]$.
By Lemma~\xref{lemma-PPP-n} $D\in \PPP_n$. Hence,
$D^+\ge D$, i.e., $D^+=D+D'$, where
$D'\ge 0$.
Let $F\subset f^{-1}(o)$ be a reduced irreducible component. Since $K_X+D$ is $\RR$-linearly trivial over
$Z$, $D'$ is vertical and $D'= c_of^*P$. Let $d_F$
and $\mu$ be multiplicities of $F$ in $D$ and $f^*o$, respectively
($\mu$ is a positive integer). Since $(X,D+D')$ is lc and $n(D+D')$ is
an integral divisor, the multiplicity of $F$ in $D+D'$
has the form $k/n$, where $k\in\ZZ$, $1\le k\le n$. Then
$k/n=d_F+c_o\mu$ and
\[
c_o=\frac1\mu\left(\frac kn-d_F\right),\qquad
d_o=1-\frac1\mu\left(\frac kn-d_F\right).
\]
Consider two cases.

a) $d_F\in \Phi(\R)$, so $d_F=1-r/m$ ($r\in\R$, $m\in\ZZ$, $m>0$).
Then we can write
\[
d_o=1-\frac{km+rn-nm}{nm\mu}=1-\frac{r'}{m\mu}<1,
\]
where
\[
0\le r'=r+\frac{km}n-m= \frac{km+rn-nm}{n}\le \frac{nm+rn-nm}n\le 1.
\]
Therefore,
$d_o\in \Phi(\R(n))$, where $0\le r'=r+\frac{km}n-m\le 1$.
This proves, in particular, (i).

b) $d_F>1-\ep _{d-1}$. In this case,
\[
1>d_o=1-\frac1\mu\left(\frac kn-d_F\right)>
1-\frac1\mu\left(\frac kn-1+\ep _{d-1}\right)>
1-\ep _{d-1}.
\]
This finishes the proof of (ii).
\end{proof}

\begin{proposition}
\label{prop-D-not-big}
Fix a finite rational subset $\R\subset [0,\, 1]$ and a positive integer 
$I$ divisible by $I(\R)$.
Let $(X,D)$ be a log semi-Fano variety of dimension $d$ such that
$X$ is $\QQ$-factorial FT and $D\in \Phi(\R)$. Assume that
there is a $(K+D)$-trivial contraction $f\colon X\to Z$ with
$0<\dim Z<d$.
Fix the choice of $I_0$ and $\psi$ in \xref{construction-adj-def-mod}
so that $\mathbf D\smod$ is effective.
We take $I$ so that $I_0$ divides $I$.
Assume the LMMP in dimension $d$.
Further, assume that Conjectures~\xref{BAB} and \xref{conj-main-adj}
hold in dimension $d-1$ and $d$, respectively.
If $K_Z+D\sdiv+D\smod$ is $Im$-complemented, then so is $K_X+D$.
\end{proposition}

\begin{proof}
Put $D_Z:=D\sdiv+D\smod$.
Apply (i) of Conjecture \xref{conj-main-adj}
to $(X,D)$. We obtain
\[
K+D=f^*(K_Z+D_Z),
\]
and $(Z,D_Z)$ is lc, where $D\sdiv\in\Phi(\R')$. 
By \eqref{conj-main-adj-2} $I''D\smod$ is integral 
for some bounded $I''$. Thus
replacing $\R'$ with $\R'\cup \{1/I'', 2/I'',\dots, (I''-1)/I''\}$
we may also assume that $D\smod \in \Phi(\R')$.
Then $D_Z \in \Phi(\R')$.
Furthermore, by Lemma \ref{lemma-FT} $Z$ is FT and by
the construction, $-(K_Z+D_Z)$ is nef.
By our inductive hypothesis $K_Z+D_Z$ has bounded complements.

Let $K_Z+D_Z^+$ be an $n$-complement of $K_Z+D_Z$ such that $I\mid n$.
 Then $D_Z^+\ge D_Z$ (see Lemmas \ref{lemma-PPP-n-1} and \ref{lemma-PPP-n}). 
Put $H_Z:=D_Z^+-D_Z$
and $D^+:=D+f^*H_Z$. Write $D^+=\sum d_i^+D_i$.
By the above, $d_i^+\ge d_i$.
We claim that $K+D^+$ is an $n$-complement of
$K+D$. Indeed, since $K+D\q{I} f^*(K_Z+D_Z)$, we have
\begin{align*}
n(K+D^+)=&n(K+D+f^*H_Z)=
\\
&(n/I)I(K+D)+(n/I)If^*H_Z\sim
\\
&(n/I)f^*I(K_Z+D_Z)+(n/I)f^*IH_Z=
\\
&(n/I)f^*I(K_Z+D_Z^+)=
f^*n(K_Z+D_Z^+)\sim f^*0=0.
\end{align*}
Thus, $n(K+D^+)\sim 0$. Further, since $nd_i^+$ is a nonnegative 
integer
and $d_i^+\ge d_i$,
the inequality
\[
nd_i^+= \down{(n+1)d_i^+}\ge \down{(n+1)d_i}
\] 
holds for every $i$ such that $0\le d_i<1$. 
Finally, by Corollary \ref{cor_Inv_aDj} 
the log divisor $K+D^+=f^*(K_Z+D_Z)$ is lc.
This proves our proposition.
\end{proof}

\subsection{Proof of Theorem \ref{main-result0}
in the case when $(X,B)$ is not klt (continued)}
\label{subsetion-0-pairs-proof}
To finish the proof Theorem \ref{main-result0}
in the non-klt case we have to consider the following situation
(see \ref{subsect-1}).
$(X',B')$ is a non-klt $0$-pair such that $B'\in \Phi(\R)$,
$X'$ is $\lambda$-lt and $X'$ is FT, where $\lambda$ 
depends only on $\R$ and the dimension of $X'$.
Moreover, there is a Fano fibration $X'\to Z'$
with $0<\dim Z'<\dim X'$.
The disired bounded $nI(\R)$-complements
exist by Proposition \ref{prop-D-not-big}
and inductive hypothesis.

\subsection{Proof of Theorem \ref{main-result} in Case \eqref{case-nef}}
\label{subsect-1-1}
To finish our proof of the main theorem we have to consider
the case when $(X,B)$ is klt and general
reduction from Section \ref{sect-reduction} leads to case \eqref{case-nef},
i.e., $-(K_Y+D_Y)$ is nef.
Replace $(X,B)$ with $(Y,B_Y)$ and put $D:=D_Y$.
Recall that in this situation $X$ is FT and
$B\in \Phi(\R, \ep')$, where $0<\ep'\le \ep_{d-1}(\ov \R)$. By 
\eqref{reduct-FT-Lambda} there is a boundary
$\Theta\ge B$ such that 
$(X,\Theta)$ is a klt $0$-pair.
For the boundary $D$ defined by \eqref{eq-D} 
we also have $D\in \Phi(\R)$ and $\down {D}\neq 0$ 
by \eqref{eq-discr-2-B-ep}. 
All these properties are preserved under birational transformations in 
\ref{subs-2-inductive-steps}.
By our assumption at the end we have case \eqref{case-nef},
i.e., $-(K+D)$ is nef (and semiample). 
Therefore it is sufficient
to prove the following.

\begin{proposition}
\label{prop-inductive}
Fix a finite rational subset $\R\subset [0,\, 1]$.
Let $(X,D=\sum d_i D_i)$ is a $d$-dimensional log semi-Fano variety
such that
\begin{enumerate}
\item
$D\in \Phi(\R)$, $(X,D)$ is not klt and $X$ is FT,
\item
there is boundary $B=\sum b_i D_i\le D$ such that
either $b_i=d_i< 1-\ep'$ or $b_i\ge 1-\ep'$ and $d_i=1$, 
where $0<\ep'\le \ep_{d-1}(\ov \R)$,
\item
$(X,\Theta)$ is a klt $0$-pair for some $\Theta\ge B$.
\end{enumerate}
Assume the LMMP in dimension $d$.
Further, assume that Conjectures~\xref{BAB} and \xref{conj-main-adj}
hold in dimension $d$.
Then $K+D$ has a bounded $n$-complement such that $I(\R) \mid n$.
\end{proposition}

The idea of the proof is to reduce the problem to 
Proposition \ref{prop-D-not-big} by considering the 
contraction $f\colon X\to Z$ 
given by $-(K+D)$. But here two technical
difficulties arise. First it may happen that the divisor
$-(K+D)$ is big and then $f$ is birational. 
In this case one can try to extend 
complements from $\down D$ but the pair $(X,D)$ is not necessarily 
plt and the inductive step (Proposition \ref{prodolj}) does not work.
We have to make some perturbations and birational 
transformations.
Second to apply inductive hypothesis to $(Z,D\sdiv+D\smod)$
we have to check if this pair satisfies conditions of 
Theorem \ref{main-result}. In particular, 
we have to check the klt property of $(Z,D\sdiv+D\smod)$.
By Corollary \ref{cor_Inv_aDj} this holds if 
any lc centre of $(X,D)$ dominates $Z$. Otherwise 
we again need some additional work.

\begin{proof}
Note that we may replace $B$ with $B_t:=tB+(1-t)D$ for $0<t<1$.
This preserves all our conditions (i)-(iii). 
Indeed, (i) and (ii) are obvious. 
For (iii), we note that $(X,D^{\lozenge})$ is a 
$0$-pair for some $D^{\lozenge}\ge D$ 
(because $-(K+D)$ is semiample).
Hence one can replace $\Theta$
with $\Theta_t:=t\Theta+(1-t)D^{\lozenge}$.

Let $\mu\colon (\tilde X,\tilde D)\to (X,D)$ be a dlt 
modification of $(X,D)$. By definition, $\mu$ is a 
$K+D$-crepant birational extraction such that
$\tilde X$ is $\QQ$-factorial, 
the pair $(\tilde X,\tilde D)$ is dlt, and each 
$\mu$-exceptional divisor $E$ has discrepancy $a(E,X,D)=-1$
(see, e.g., \cite[21.6.1]{Utah}, \cite[3.1.3]{Prokhorov-2001}).
In particular, $\tilde D\in \Phi(\R)$ and $\tilde X$ is FT by 
Lemma \ref{lemma-FT}. Let $\tilde B$ be the crepant pull-back of $B$.
One can take $t$ so that the multiplicities in $\tilde B$ 
of $\mu$-exceptional divisors are $\ge 1-\ep_{d-1}$.
Thus for the pair $(\tilde X,\tilde D)$ conditions (i)-(iii) 
hold. Therefore, we may
replace $(X,D)$ with $(\tilde X,\tilde D)$
(and $B$, $\Theta$ with their crepant pull-backs).

Let $f\colon X\to Z$ be the contraction given by $-(K+D)$.
By Theorem \ref{main-result0} (see \ref{subsect-1} and \ref{subsetion-0-pairs-proof}) 
we may assume that $\dim Z>0$.
We apply induction by $N:=\dim X-\dim Z$.

First, consider the case $N=0$. Then $-(K+D)$ is big.
We will show that $K+D$ is $n$-complemented for
some $n\in \NNN_{d-1}(\ov \R)$.

Fix $n_0\gg 0$, and let $\delta:=1/n_0$. 
Then $D_\delta:=D-\delta \down D\in \Phi(\R)$.
It is sufficient to show that $K_X+D_\delta$ is $n$-complemented for
some $n\in \NNN_{d-1}(\ov \R)$.
We will apply a variant of \cite[Th. 5.1]{Prokhorov-Shokurov-2001}
with hyperstandard multiplicities.
To do this, we run $-(K+D_\delta)$-MMP over $Z$.
Clearly, this is equivalent $\down D$-MMP over $Z$.
This process preserve the $\QQ$-factoriality and 
lc (but not dlt) property of $K+D$.
At the end we get a model $(X',D')$ such that $-(K_{X'}+D'_\delta)$
is nef over $Z$. Since $X'$ is FT, the Mori cone $\NE (X')$ is rational
polyhedral. Taking our condition $0<\delta\ll 1$ into account we get that
$-(K_{X'}+D'_\delta)$ is nef. Since
\[
-(K_{X'}+D'_\delta)= -(K_{X'}+D')+\delta \down {D'},
\]
where $\down {D'}$ is effective, $-(K_{X'}+D'_\delta)$ is also big.
Note that $(X',D')$ is lc but not klt.
By our assumptions,
\[
D'_\delta=(1-\delta)D'+\delta B' \le (1-\delta)D'+\delta\Theta'
\]
and $(X',(1-\delta)D'+\delta\Theta')$ is klt.
Therefore so is $(X',D'_\delta)$.
Now we apply \cite[Th. 5.1]{Prokhorov-Shokurov-2001}
with $\Phi=\Phi(\R)$. 
This says that we can extend complements from 
some (possibly exceptional) divisor. By 
Proposition \ref{prop-prop} the multiplicities 
of the corresponding different are contained in 
$\ov \R$.
We obtain an $n$-complement of $K_{X'}+D'_\delta$ for
some $n\in \NNN_{d-1}(\ov \R)$. By Proposition \ref{cor-pull-back_compl-I}
we can pull-back this complement to $X$ (we use the inclusion
$D'_\delta\in \Phi(\R)\subset \PPP_n$).

Now assume that Proposition \ref{prop-inductive} holds for all $N'<N$.
Run $\down D$-MMP over $Z$. After some flips and divisorial contractions
we get a model on which $\down D$ is nef over $Z$.
Since $X$ is FT, the Mori cone $\NE(X)$ is rational polyhedral.
Hence $-(K+D-\delta \down D)$ is nef for $0<\delta \ll 1$.
As above, put $D_\delta :=D-\delta \down D$. We can take $\delta =1/n_0$, $n_0\gg 0$
and then $D_\delta \in \Phi(\R)$.
On the other hand, $D_\delta\le (1-\delta)D +\delta B$ for some $\delta>0$.
Therefore, $(X,D_\delta)$ is klt. Now let $f^\flat\colon X\to Z^\flat$
be the contraction given by $-(K+D_\delta)$. Since
$-(K+D_\delta)=-(K+D)+\delta \down D$, there is decomposition
$f\colon X \stackrel {f^\flat} \longrightarrow Z^\flat \longrightarrow Z$.

If $\dim Z^\flat=0$, then $Z^\flat=Z$ is a point, a contradiction.
If $\dim Z^\flat<\dim X$, then by Corollary \ref{cor_Inv_aDj}
$(Z^\flat, (D_\delta)\sdiv+(D_\delta)\smod)$ is a klt
log semi-Fano variety. We can apply Proposition \ref{prop-D-not-big}
to the contraction $X\to Z^\flat$ 
and obtain a bounded complement of $K+D_\delta$.
Clearly, this will be a complement of $K+D$.

Therefore, we may assume that $-(K+D_\delta)$ is big, $f^\flat$ is birational,
and so $\down D$ is big over $Z$.
In particular, the horizontal part $\down{D}^{\h}$
of $\down D$ in non-trivial.

Replace $(X,D)$ with its dlt modification.
Assume that $\down{D}^{\h}\neq\down D$.
As above, run $\down{D}^{\h}$-MMP over $Z$. For $0<\delta \ll 1$,
the divisor $-(K+D-\delta \down{D}^{\h})$ will be nef. Moreover,
it is big over $Z$. Therefore, $-(K+D-\delta \down{D}^{\h})$ defines
a contraction $f'\colon X\to Z'$ with $\dim Z'>\dim Z$.
By our inductive hypothesis there is a bounded complement.

It remains to consider the case when $\down{D}^{\h}=\down D$.
Then any lc centre of $(X,D)$ dominates $Z$.
By Corollary \ref{cor_Inv_aDj} and Proposition \ref{prop-D-not-big}
there is a bounded complement of $K+D$.

\end{proof}
This finishes the proof of Theorem \ref{main-result}.
Corollaries \ref{main-result-3-Cor-1-b} and \ref{main-result-2}
immediately follows by this theorem, Corollary \ref{corollary-815},
and \cite{Alexeev-1994}.

\subsection{Proof of Corollary \ref{main-result-3-Cor-1-aC}}
\label{sketch-proof-main-result-3-Cor-1-aC}
Replacing $(X,D)$ with its log terminal modification
we may assume that $(X,D)$ is dlt.
If $D=0$, 
we have $nK\sim 0$ for some $n\le 21$
by \cite{Blache-1995}. Thus we assume that 
$D\neq 0$. Run $K$-MMP. 
We can pull-back complements by 
Proposition \ref{cor-pull-back_compl-I}.
The end result is a $K$-negative extremal contraction $(X',D')\to Z$
with $\dim Z\le 1$.
If $Z$ is a curve, then either $Z\simeq\PP^1$ or $Z$ is an elliptic curve.
In both cases we apply Proposition \ref{prop-D-not-big}
(the FT property of $X$ is not needed).
Otherwise $X'$ is a klt log del Pezzo surface with $\rho(X')=1$.
In particular, $X'$ is FT. In this case the assertion follows by 
Theorem \ref{main-result0}.

\subsection{Proof of Corollary \ref{main-result-3-Cor-1-a}}
\label{sketch-proof-main-result-3-Cor-1-a}
First we construct a 
crepant dlt model $(\bar X,\bar D)$ of $(X,D)$ such that 
each component of $\bar D$ meets $\down {\bar D}$.
Replacing $(X,D)$ with its log terminal modification
we may assume that $(X,D)$ is dlt, $X$ is $\QQ$-factorial, and
$\down D\neq 0$. 
If $\down {D}=D$, we put $(\bar X,\bar D)=(X,D)$.
Otherwise, run $K+D-\down D$-MMP.
Note that none of connected components of $\down D$ 
is contracted. Moreover, the number of connected components 
of $\down D$ remains the same (cf. \cite[Prop. 12.3.2]{Utah},
\cite[Th. 6.9]{Shokurov-1992-e}). 
At the end we get an extremal contraction $(X',D')\to Z$
with $\dim Z\le 2$.
If $Z$ is not a point, we can apply Proposition \ref{prop-D-not-big}.
Otherwise, $\rho(X')=1$, $\down{D'}$ is connected, and 
each component of $D'$ meets $\down{D'}$. 
The same holds on a log terminal modification
$(\bar X,\bar D)$ of $(X',D')$
because $X$ is $\QQ$-factorial.

By Corollary \ref{main-result-3-Cor-1-aC},
for each component $\bar D_i\subset \down{\bar D}$,
the log divisor $K_{\bar D_i}+\Diff_{\bar D_i}(\bar D-\bar D_i)$
has bounded complements, i.e., 
there is $n_0=n_0(\R)$ such that
$n_0(K_{\bar D_i}+\Diff_{\bar D_i}(\bar D-\bar D_i))\sim 0$.
Thus we may assume that $n_0(K_{\bar X}+\bar D)|_{\down{\bar D}}\sim 0$.
Recall that the multiplicities of $\Diff_{\bar D_i}(\bar D-\bar D_i)
=\sum \delta_j\Delta_j$
are computed by the formula $\delta_j=1-1/m_j+(\sum_l k_l d_l)/m_j$,
where $m_j,\, k_l\in \ZZ$, $m_j>0$, $k_l\ge 0$, 
$d_l$ are multiplicities of $\bar D$, and $\sum_l k_l d_l\le 1$.
Since $d_l\in \Phi(\R)$, there is only a finite number of 
possibilities for $d_l$ with $k_l\neq 0$.
Further, since $n_0\delta_j\in \ZZ$, there is only a finite number of 
possibilities for $m_j$. Thus we can take $n_1=n_1(\R)$ such that
$n_1\bar D$ is an integral divisor and
$n_1(K_{\bar X}+\bar D)|_{\down{\bar D}}\sim 0$.
Since $\bar X$ is FT, there is an integer $n_2$ such that 
$n_1n_2(K_{\bar X}+\bar D)\sim 0$ on $\bar X$.
This defines a cyclic \'etale over $\down{\bar D}$ cover 
$\pi\colon \hat X\to \bar X$.
Let $\hat D:=\pi^*\bar D$. Then 
$(\hat X,\hat D)$ is a $0$-pair such that
$\down {\hat D}$ has at least $n_2$ connected components.
On the other hand, the number of connected componets 
of a $0$-pair is at most two (see
\cite[2.1]{Fujino-2000-ab}, cf. \cite[6.9]{Shokurov-1992-e}).
Thus, $n(K_{\bar X}+\bar D)\sim 0$, where $n=2n_1$.
This proves our corollary.

\subsection{Proof of Corollary \ref{main-result-3-Cor-1-2div}}
\label{sketch-proof-main-result-3-Cor-2div}
In notation of \ref{eq-reduction-setup},
take $0<\bar \ep<\ep_{2}(\ov \R)/2$.
We may assume that $(X,B)$ is such as in \ref{sub-reduction-1},
so there are (at least) two components $B_1$ and $B_2$ of multiplicities 
$b_i\ge 1-\bar \ep$ in $B$. 
Then by Lemma \ref{pair-discr} and
\eqref{eq-discr-2-B-ep}
components $B_1$, $B_2$ do not meet each other
and by Corollary \ref{cor-reduct-discrs} this holds on each step
of the LMMP as in \ref{subs-2-inductive-steps}. 
Therefore, we cannot
get a model with $\rho=1$. In particular, case
\eqref{case-rho=1} is impossible.

Consider case \eqref{case-nef}.
If the divisor $-(K_Y+D_Y)$ is big, we can argue 
as in the proof of Proposition \ref{prop-inductive}.
Then we do not need Conjecture \ref{BAB}.
If $K_Y+D_Y\equiv 0$, we
can use Corollary \ref{main-result-3-Cor-1-a}; it is sufficient 
to have only
one divisor $E$ (exceptional or not) with $a(E,X,D)\le -1+\bar \ep$.
In other cases we use induction 
to actual fibrations (Proposition \ref{prop-D-not-big}), that is,
with the fibres and the base of dimension $\ge 1$ and
by our assumptions with dimensions $\le 2$.


\providecommand{\bysame}{\leavevmode\hbox to3em{\hrulefill}\thinspace}
\providecommand{\MR}{\relax\ifhmode\unskip\space\fi MR }
\providecommand{\MRhref}[2]{%
  \href{http://www.ams.org/mathscinet-getitem?mr=#1}{#2}
}
\providecommand{\href}[2]{#2}

\end{document}